\documentclass[11pt]{article}
\usepackage{amssymb}
\usepackage{amsthm}

\usepackage{verbatim}

\usepackage{mathrsfs}
\usepackage{graphicx,subfigure}
\usepackage{paralist} 
\usepackage{enumerate}
\usepackage{hyperref,url}

\usepackage{bm} 
\usepackage{placeins}
\usepackage{color}
\usepackage[margin=30mm]{geometry}

\newcommand{\pd}{\partial}
\newcommand{\abs}[1]{\left| #1 \right|}

\newcommand{\vel}{\bm{v}}
\newcommand{\bbv}{\bm{v}}
\newcommand{\vtau}{\bm{v_\tau}}

\newcommand{\eps}{\ensuremath{\varepsilon}}

\newcommand{\mgp}{\ensuremath{| \nabla\phi |}}
\newcommand{\mgpp}{\ensuremath{|\nabla\phi^{m}|}}
\newcommand{\mgpn}{\ensuremath{| \nabla\phi^{m-1}|}}

\newcommand{\io}{\ensuremath{\int_\Omega}}

\newcommand{\up}{\ensuremath{u^m_h}}
\newcommand{\un}{\ensuremath{u^{m-1}_h}}

\newcommand{\fe}{\ensuremath{f^e}}
\newcommand{\ue}{\ensuremath{u^e}}
\newcommand{\ep}{\ensuremath{e_h^{m}}}
\newcommand{\en}{\ensuremath{e_h^{m-1}}}

\newcommand{\nabphi}{\ensuremath{\nabla_{\phi}}}

\newcommand{\rhop}{\ensuremath{\rho^m}}
\newcommand{\rhon}{\ensuremath{\rho^{m-1}}}

\newcommand{\vp}{\ensuremath{\bm{v}^{m}}}
\newcommand{\tint}{\ensuremath{\int_{t_{m-1}}^{t_m}}}

\newcommand{\ove}{\ensuremath{\frac1{\varepsilon}}}
\newcommand{\ldot}{.}
\newcommand{\dhp}{\ensuremath{d^{m}}}
\newcommand{\dhn}{\ensuremath{d^{m-1}}}
\newcommand{\be}{\begin{equation}}
\newcommand{\ee}{\end{equation}}
\newcommand{\been}{\begin{eqnarray*}}
\newcommand{\eeen}{\end{eqnarray*}}
\newcommand{\nono}{\nonumber}

\theoremstyle{plain}
\newtheorem{theorem}{Theorem}[section]

\newtheorem{lemma}[theorem]{Lemma}
\newtheorem{corollary}[theorem]{Corollary}
\newtheorem{remark}[theorem]{Remark}

\theoremstyle{plain}

\title{Stability and error analysis for a diffuse interface approach to an advection–-diffusion equation on a moving surface}

\author{Klaus Deckelnick  \and Vanessa Styles}

\author{Klaus Deckelnick\footnote{Institut f\"ur Analysis und Numerik,
Otto-von-Guericke-Universit\"at Magdeburg, Universit\"atsplatz 2, 39106 Magdeburg, Germany.}
 ~and Vanessa Styles\footnote{Department of Mathematics,
University of Sussex, Brighton BN1 9QH, UK.} }

\date{~}

\begin{document}

\maketitle


\begin{abstract}
In this paper we analyze a fully discrete numerical scheme for solving a parabolic PDE on a moving surface. The method
is based on a diffuse interface approach that involves a level set description of the moving surface. Under suitable
conditions on the spatial grid size, the time step and the interface width we obtain stability and error bounds
with respect to natural norms. Furthermore, we present test calculations that confirm our analysis.
\end{abstract}

\section{Introduction}
Let $\{ \Gamma(t) \}_{t \in [0,T]} $ be a family of closed hypersurfaces in $\mathbb{R}^{n+1} (n=1,2)$ evolving in time. 
In this paper we consider a finite element approach for solving the advection–-diffusion
equation 
\begin{eqnarray} 
\partial_t^{\bullet}  u  + u \nabla_{\Gamma} \cdot \vel -  \Delta_{\Gamma} u  & = & f \qquad \mbox{on } S_T \label{surfeq} \\
u(\cdot,0) & = & u_0 \quad \; \mbox{ on } \Gamma(0), \label{surfinit}
\end{eqnarray}
where $S_T= \bigcup_{t \in (0,T)} \bigl( \Gamma(t) \times \lbrace t \rbrace \bigr)$  and $\vel: \overline{S_T} \rightarrow \mathbb{R}^{n+1}$  denotes a given velocity field.
Furthermore, $\nabla_{\Gamma}$ is the tangential gradient, $\Delta_{\Gamma} = \nabla_{\Gamma} \cdot \nabla_{\Gamma}$ the Laplace Beltrami operator and
$\partial_t^{\bullet}=\partial_t + \vel \cdot \nabla$ denotes the material derivative. \\
Parabolic surface PDEs of the form (\ref{surfeq}) have applications in fluid dynamics and materials science, such as the transport and diffusion of surfactants on a fluid/fluid interface, 
\cite{S90} or
diffusion-induced grain boundary motion, \cite{CFP97}. In these as in several other applications the velocity $\vel$ is not given but determined through an additional equation so that
(\ref{surfeq}) becomes a subproblem of a more complicated system in which the variable $u$ is coupled to other variables. The analysis and the numerical solution of such systems
then naturally requires the development of corresponding methods for (\ref{surfeq}). We refer to \cite{DzE13} for a comprehensive overview of finite element methods for
solving PDEs on stationary and evolving surfaces. \\
Concerning the numerical methods that have been proposed for (\ref{surfeq}) one may distinguish between Lagrangian and
Eulerian type schemes. The first approach has been pursued by Dziuk and Elliott within their evolving surface finite element method,  \cite{DE07}, which uses polyhedral approximations of the evolving hypersurfaces
$\Gamma(t)$. While \cite{DE07} contains an error analysis in the  spatially discrete case, the fully discrete case is investigated in \cite{DE12}, \cite{DLM12} and \cite{LMV13}. Optimal $L^2$-error bounds are 
obtained in \cite{DE13} and a corresponding finite volume approach is proposed and analyzed in \cite{LNR11}. Since the mesh for the discretization of (\ref{surfeq}) is fitted to the hypersurface
$\Gamma(t)$, a coupling to a bulk equation is not straightforward. This difficulty is not present in Eulerian type schemes, in which $\Gamma(t)$ is typically described via a level set function
defined in an open neighbourhood of  $\Gamma(t)$. In order to discretize the surface PDE in this setting it has been
proposed in \cite{AS03}, \cite{BCOS} and \cite{XZ03} to extend the surface quantity $u$ to a band around $\Gamma(t)$ and to solve a suitable (weakly) parabolic PDE in that bulk region using a finite 
difference method. In \cite{DE08} and \cite{DE10}, the same idea is used in a finite element context for which the underlying variational formulation is derived with the help of
a transport identity. An Eulerian finite element approach that doesn't use an extended PDE is proposed and analyzed in \cite{ORX14} and \cite{OR14}. The method is based on a weak formulation on the
space-time manifold and the finite element space is obtained by taking  traces of the
corresponding bulk finite elements. The approximation of $\Gamma(t)$ on which these spaces are defined usually  arises from a suitable interpolation 
of the given level set function describing $\Gamma(t)$.  The resulting discrete hypersurface will in general cut arbitrarily through the background mesh and its location forms one of the
main difficulties in implementing the scheme. A different approach of generating the discrete hypersurfaces is pursued in \cite{HLZ15}, where a discretization of  (\ref{flowmap}) below is 
combined with the cut finite element technique. Finally, Section 5 in \cite{DER14} proposes a hybrid method that employs the above--mentioned idea of trace finite elements
together with a narrow band technique for the elliptic part of the PDE. \\
In this paper we are concerned with the diffuse interface approach for solving (\ref{surfeq}), which was introduced in \cite{RV06} for a stationary surface and
in \cite{ESSW11} and \cite{TLLWV09} for evolving surfaces. As in some of the methods described above, 
the surface quantity $u$ is extended to  a bulk quantity satisfying a suitable parabolic PDE in a neighbourhood of $\Gamma(t)$ and
the bulk equation is then localized to a thin layer
of thickness $\epsilon$ with the help of a phase field function (see \cite{ES09} for a corresponding convergence analysis). Since we are
interested in using finite elements, the localized PDE needs to be written in a suitable variational form. Following \cite{ESSW11} this is achieved 
with the help of a transport identity and results in a discretization by linear finite elements in space and a backward Euler scheme in time. The 
detailed derivation along with an existence result for the discrete solution will be given in Section 3.
An advantage of this approach is that in the implementation the evolution of the hypersurfaces is easily
incorporated by evaluating the phase field function. As the main new contribution of our paper
 we shall derive conditions relating
the interface width $\epsilon$, the spatial grid size $h$ and the time step $\tau$ which allow for a rigorous stability and error analysis.
The corresponding results are formulated and proved in Sections 4 and 5 respectively, while we report in Section 6 on results of numerical tests. 
Let us finally remark that a phase field approach involving a nonlocal phase field function and finite elements has been proposed in \cite{Bu09} for an elliptic
surface PDE. Theorem 7 in \cite{Bu09} provides an error estimate in terms of an approximation error and an error due to the
phase field representation. The latter decays at a rate $O(\epsilon^p)$ for some $p<1$, while a coupling between $\epsilon$ and 
the grid size $h$ is not discussed.

\section{Preliminaries}
\subsection{Surface representation and surface derivatives}
For each $t \in [0,T]$ let  $\Gamma(t) \subset \mathbb{R}^{n+1} \, (n=1,2)$ be a connected, compact and orientable hypersurface without boundary. We
suppose that $\vel: \overline{S_T} \rightarrow \mathbb{R}^{n+1}$ is a prescribed velocity field of the form
\begin{equation}  \label{velsplit}
\vel = V \nu + \vtau, \quad \mbox{ with } (\vtau, \nu)=0.
\end{equation}
Here,  $\nu$ is a unit normal and $V$ the corresponding normal velocity of $\Gamma(t)$ and $(\cdot, \cdot)$ denotes the Euclidian scalar product in $\mathbb{R}^{n+1}$. 
Note that the normal part  $V\nu$ is responsible for the geometric motion of $\Gamma(t)$, while the tangential part $\vtau$ is associated with the transport of material along the surface.
We assume that there exists a smooth map 
$\Phi:\Gamma(0) \times [0,T] \rightarrow \mathbb{R}^{n+1}$ such that $\Phi(\cdot,t)$ is a diffeomorphism from $\Gamma(0)$ onto
$\Gamma(t)$ for every $t \in [0,T]$ satisfying
\begin{eqnarray}  
\frac{\partial \Phi}{\partial t}(P,t) & = &  \vel(\Phi(P,t),t), \quad P \in \Gamma(0), t \in (0,T]; \label{flowmap} \\
\Phi(P,0) & = & P, \qquad \qquad \quad \; \,  P \in \Gamma(0).  \label{flowinit}
\end{eqnarray}
Let us next introduce the differential operators which are required to formulate our PDE. To begin, for fixed $t$ and a function $\eta: \Gamma(t) \rightarrow \mathbb{R}$ we denote
by $\nabla_{\Gamma} \eta=(\underline{D}_1 \eta,\ldots, \underline{D}_{n+1}\eta)$ its tangential gradient. If $\bar{\eta}$ is an extension of $\eta$ to an open neighbourhood
of $\Gamma(t)$ then
\begin{equation}  \label{tanggrad0}
\nabla_{\Gamma} \eta(x) = \bigl( I- \nu(x,t) \otimes \nu(x,t) \bigr) \nabla \bar{\eta}(x), \quad x \in \Gamma(t).
\end{equation}
Furthermore, $\Delta_{\Gamma} \eta= \nabla_{\Gamma} \cdot \nabla_{\Gamma} \eta = \sum_{i=1}^{n+1} \underline{D}_i \underline{D}_i \eta$ denotes the Laplace–-Beltrami operator. \\[2mm]
Next, for a smooth function $\eta$ on $S_T$ we define the material derivative of $\eta$ at $(x,t)=(\Phi(P,t),t)$ by 
$\partial_t^{\bullet}\eta(x,t):= \frac{d}{dt} [ \eta(\Phi(P,t),t)]$. If $\bar{\eta}$ is an extension of $\eta$ to an open space-time neighbourhood, then
\begin{displaymath}
\partial_t^{\bullet} \eta(x,t) = \bar{\eta}_t(x,t) + (\vel(x,t),\nabla \bar{\eta}(x,t)), \quad (x,t) \in S_T.
\end{displaymath}

\vspace{2mm}

Our numerical approach will be based on an implicit representation of $\Gamma(t)$, so that we suppose in what follows that there exists a 
smooth function $\phi: \Omega \times [0,T] \rightarrow \mathbb{R}$ such that for $0 \leq t \leq T$
\begin{equation}  \label{implicit}
\displaystyle
\Gamma(t) = \lbrace x \in \Omega \, | \, \phi(x,t) = 0 \rbrace \quad \mbox{ and } \quad \nabla \phi(x,t) \neq 0, \, x \in \Gamma(t).
\end{equation}
Here, $\Omega \subset \mathbb{R}^{n+1}$ is a bounded domain with $\Gamma(t) \subset \Omega$ for all $t \in [0,T]$. For later use we introduce for $t \in [0,T], \,
r>0$ the sets
\begin{displaymath}
U_r(t):= \lbrace x \in \Omega \, | \, | \phi(x,t) | <r \rbrace \quad \mbox{ and } \quad \mathcal U_{r,T}:= \bigcup_{t \in [0,T]} \bigl( U_r(t) \times \lbrace t \rbrace \bigr).
\end{displaymath}
In view of (\ref{implicit}) there exist $\delta_0>0, 0< c_0 \leq c_1, \, c_2>0$ such that $\overline{U_{\delta_0}(t)} \subset \Omega, 0 \leq t \leq T$ and 
\begin{equation} \label{gradbound}
\displaystyle c_0 \leq | \nabla \phi(x,t) | \leq c_1, \;   | D^2 \phi(x,t) |, \, | \phi_t(x,t) |, \, | \phi_{tt}(x,t) | \leq c_2,  \quad (x,t)  \in  \mathcal U_{\delta_0,T}.
\end{equation}

\subsection{Extension} \label{extension}
Our next aim is to extend functions defined on $S_T$ to a space-time neighbourhood. A common approach which is
well suited to a description of $\Gamma(t)$ via the signed distance function consists in extending constantly in normal direction. In what follows we shall introduce
a suitable generalization to the case (\ref{implicit}). 
Consider for $P \in \Gamma(0)$ and $t \in [0,T]$ the parameter-dependent system of ODEs
\begin{equation}
  \label{odesys}
  \gamma_{P,t}'(s) = \frac{\nabla \phi(\gamma_{P,t}(s),t)}{| \nabla \phi(\gamma_{P,t}(s),t) |^2}, \quad \gamma_{P,t}(0)= \Phi(P,t).
\end{equation}
Using a compactness argument it can be shown that there exists $0< \delta < \delta_0$ so that the solution $\gamma_{P,t}$ of (\ref{odesys}) exists uniquely on $(-\delta,\delta)$ 
uniformly in $P \in \Gamma(0), t \in [0,T]$. 
Thus we can define the smooth mapping $F_t:\Gamma(0) \times (-\delta,\delta) \rightarrow \mathbb{R}^{n+1}$ by
\begin{equation}
  \label{defF}
  F_t(P,s):=\gamma_{P,t}(s), \quad P \in \Gamma(0), |s| < \delta.
\end{equation}
In view of the chain rule and (\ref{odesys}) we immediately see that  $\frac{d}{ds} \phi(\gamma_{P,t}(s),t) = 1$, which implies that $\phi(\gamma_{P,t}(s),t)=s, |s| < \delta$ since
$\gamma_{P,t}(0)=\Phi(P,t) \in \Gamma(t)$. In particular,  $x=F_t(P,s)$ yields that $| \phi(x,t) | < \delta$ and it is not difficult to verify  that $F_t$ is a diffeomorphism of $\Gamma(0) \times (-\delta,\delta)$ 
onto $U_{\delta}(t)$ for $t \in [0,T]$, whose inverse has the form
\begin{equation}
  \label{eq:Finv}
  F_t^{-1}(x)= (p(x,t),\phi(x,t)), \quad x \in U_{\delta}(t).
\end{equation}
Here,  $p: \mathcal U_{\delta,T} \rightarrow \mathbb{R}^{n+1}$ satisfies $p(x,t) \in \Gamma(0), x \in U_{\delta}(t)$. 
Furthermore, since $\phi(F_t(P,s),t)=s$ we deduce from (\ref{eq:Finv}) that
\begin{equation}  \label{pprop}
\displaystyle p(x,t) = P, \quad \mbox{ if } x= F_t(P,s) \in U_{\delta}(t).
\end{equation}
The function 
$\tilde{p}: \mathcal U_{\delta,T} \rightarrow \mathbb{R}^{n+1}, \, \tilde{p}(x,t):= \Phi(p(x,t),t)$
then is smooth and satisfies $\tilde{p}(x,t) \in \Gamma(t), 0 \leq t \leq T$. In addition we claim that
\begin{equation} \label{tildep}
\displaystyle \tilde{p}(x,t) = x, \quad x \in \Gamma(t).
\end{equation}
To see this, let $x \in \Gamma(t)$, say $x=\Phi(P,t)=\gamma_{P,t}(0)=F_t(P,0)$ for some $P \in \Gamma(0)$.  Using (\ref{pprop}) with $s=0$ we
deduce that
\begin{displaymath}
\tilde{p}(x,t)= \Phi(p(x,t),t) = \Phi(P,t) = x,
\end{displaymath}
proving (\ref{tildep}). 
Let us next use $\tilde{p}$ in order to extend a function  
$z: \overline{S_T} \rightarrow \mathbb{R}$ to $\mathcal U_{\delta,T}$ by setting
\begin{equation} \label{ext}
\displaystyle 
z^e(x,t):= z(\tilde{p}(x,t),t), \quad (x,t) \in \mathcal U_{\delta,T}.
\end{equation}
Clearly, $z^e(\cdot,t)=z(\cdot,t)$ on $\Gamma(t)$ by (\ref{tildep}). Moreover, (\ref{pprop}) implies for $P \in \Gamma(0), |s| < \delta$
\begin{displaymath}
z^e(F_t(P,s),t) = z \bigl( \tilde{p}(F_t(P,s),t),t \bigr)=  z \bigl(\Phi(p(F_t(P,s),t),t),t \bigr) = z(\Phi(P,t),t),
\end{displaymath}
from which we obtain by differentiating with respect to $s$ and using (\ref{odesys}), (\ref{defF}) that
\begin{equation}  \label{zetang}
\displaystyle 
  \bigl( \nabla z^e(x,t),\nabla \phi(x,t) \bigr) = 0, \quad (x,t)  \in \mathcal U_{\delta,T}.
\end{equation}

\begin{lemma} \label{extest} Let $z^e$ be defined by (\ref{ext}). Then we have for $t \in [0,T], \, 0 < r < \delta$ and $| \alpha | =k \in \lbrace 0,1,2 \rbrace$:
\begin{eqnarray}
\Vert D_x^{\alpha} z^e(\cdot,t)  \Vert_{L^2(U_r(t))} &  \leq &  C \sqrt{r} \Vert z(\cdot,t) \Vert_{H^k(\Gamma(t))};   \label{zeest1} \\[2mm]
\Vert D_x^{\alpha} z^e_t(\cdot,t) \Vert_{L^2(U_r(t))} & \leq & C \sqrt{r} \bigl( \Vert \partial^{\bullet}_t z (\cdot,t) \Vert_{H^k(\Gamma(t))} + \Vert z(\cdot,t) \Vert_{H^{k+1}(\Gamma(t))} \bigr).\label{zeest2}
\end{eqnarray}
\end{lemma}
{\bf Proof:} Let us recall that $F_t$ is a diffeomorphism from $\Gamma(0)  \times (-r,r)$ onto $U_r(t)$ while $\Phi(\cdot,t)$ is a diffeomorphism from $\Gamma(0)$ onto
$\Gamma(t)$. We deduce from (\ref{pprop}) and the definition of $\tilde{p}$  that  $\tilde{p}(F_t(P,s),t)=\Phi(P,t), P \in \Gamma(0), | s| < r$ so that we obtain  with the help of
the transformation rule
\begin{equation}  \label{intest}
\displaystyle 
\int_{U_r(t)} | z^e(x,t) |^2 dx = 
\int_{U_r(t)}  |  z(\tilde{p}(x,t),t) |^2 dx \leq c \int_{-r}^r 
\int_{\Gamma(0)}  | z(\Phi(P,t),t) |^2 do_P ds \leq c r \int_{\Gamma(t)}  | z(Q,t) |^2 do_Q
\end{equation} 
which is (\ref{zeest1}) for $k=0$. Next, differentiating the identity $\phi(\tilde{p}(x,t),t)=0$ with respect to $x_i$ we infer that
$(\nabla \phi(\tilde{p}(x,t),t),\tilde{p}_{x_i}(x,t))=0,i=1,\ldots,n+1$. Hence we obtain from (\ref{ext}) and (\ref{tanggrad0}) that
\begin{eqnarray} 
z^e_{x_i}(x,t) & = & \sum_{k=1}^{n+1} z^e_{x_k}(\tilde{p}(x,t),t) \tilde{p}_{k,x_i}(x,t) = \sum_{k=1}^{n+1} \underline{D}_k z(\tilde{p}(x,t),t) \tilde{p}_{k,x_i}(x,t), \label{zederiv1} \\
z^e_{x_i x_j}(x,t) & = & \sum_{k,l=1}^{n+1} \underline{D}_l \underline{D}_k z(\tilde{p}(x,t),t) \tilde{p}_{k,x_i}(x,t) \tilde{p}_{l,x_j}(x,t) 
 + \sum_{k=1}^{n+1} \underline{D}_k z(\tilde{p}(x,t),t) \tilde{p}_{k,x_i x_j}(x,t). \label{zederiv2} 
\end{eqnarray}
Similarly, $(\nabla \phi(\tilde{p}(x,t),t),\tilde{p}_t(x,t)) = - \phi_t(\tilde{p}(x,t),t)=(\nabla \phi(\tilde{p}(x,t),t), \vel(\tilde{p}(x,t),t))$ by (\ref{velext}) below, so that
\begin{eqnarray}
z^e_t(x,t) & = &  z^e_t(\tilde{p}(x,t),t) + (\nabla z^e(\tilde{p}(x,t),t), \tilde{p}_t(x,t)) \nonumber \\
& = & \partial^{\bullet}_t z(\tilde{p}(x,t),t) + \sum_{k=1}^{n+1} \underline{D}_k z(\tilde{p}(x,t),t) \bigl( \tilde{p}_{k,t}(x,t) - \vel_k(\tilde{p}(x,t),t) \bigr).  \label{zemat}
\end{eqnarray}
Combining (\ref{zederiv1}), (\ref{zederiv2}) with the argument in (\ref{intest}) we obtain (\ref{zeest1}). The estimate (\ref{zeest2}) follows in a similar way
if one starts from (\ref{zemat}). \qed \\[2mm]

Let us next extend the surface differential operators $\nabla_{\Gamma}$ and $\partial^{\bullet}_t$. By reversing the orientation of $\Gamma(t)$ if necessary we
may assume that the functions  $\nu: \mathcal U_{\delta,T} \rightarrow
\mathbb{R}^{n+1}$, $V: \mathcal U_{\delta,T} \rightarrow \mathbb{R}$ defined by
\begin{displaymath}
\nu(x,t):= \frac{ \nabla \phi(x,t)}{| \nabla \phi(x,t) |}, \quad V(x,t):= - \frac{\phi_t(x,t)}{| \nabla \phi(x,t) |}, \quad (x,t) \in  \mathcal U_{\delta,T}
\end{displaymath}
are extensions of the unit normal and the normal velocity respectively. In particular, we define for a function $\eta \in C^1(U_{\delta}(t))$ its Eulerian tangential gradient by
\begin{equation} \label{tanggrad}
\displaystyle 
\nabla_{\phi} \eta(x):= \bigl( I - \nu(x,t) \otimes \nu(x,t) \bigr) \nabla \eta(x), \quad x \in U_{\delta}(t)
\end{equation}
and remark that $(\nabla_{\phi} \eta)_{| \Gamma(t)} = \nabla_{\Gamma} [\eta_{| \Gamma(t)}]$. Furthermore, it follows from Lemma 2 in \cite{DE10} that for $\eta \in C^1_0(\Omega)$
with $\mbox{supp}\eta \subset U_{\delta}(t)$ 
\begin{equation} \label{Eulintparts}
\displaystyle \int_{\Omega} \nabla_{\phi} \eta \, | \nabla \phi | = - \int_{\Omega} \eta H \nu \, | \nabla \phi |, \quad \mbox{ where } H= - \nabla \cdot \nu.
\end{equation}
Note that $H_{| \Gamma(t)}$ is the mean curvature of $\Gamma(t)$. \\[2mm]
Let us also extend the  velocity field $\vel$ to $\mathcal U_{\delta,T}$. We first extend
its tangential part by setting
\begin{displaymath}
{\bf \tilde{v}_{\tau}}(x,t):= (I - \nu(x,t) \otimes \nu(x,t)) \vtau^e(x,t), \quad  (x,t)  \in \mathcal U_{\delta,T}.
\end{displaymath}
In view of (\ref{velsplit}) the function $\vel(x,t):= V(x,t) \nu(x,t) + {\bf \tilde{v}_{\tau}}(x,t)$ extends the given velocity field from $\overline{S_T}$ to $\mathcal U_{\delta,T}$
and satisfies
\begin{equation} \label{velext}
\displaystyle
\phi_t + ( \vel, \nabla \phi) =0 \quad \mbox{ in } \mathcal U_{\delta,T}.
\end{equation}
In particular, we can use the exended velocity $\vel$ to define the material derivative for a function $\eta$ on $\mathcal U_{\delta,T}$
by setting
\begin{displaymath}
\partial_t^{\bullet} \eta(x,t):= \eta_t(x,t) + (\vel(x,t),\nabla \eta(x,t)), \quad (x,t)  \in \mathcal U_{\delta,T}.
\end{displaymath}

\section{Weak formulation and numerical scheme}
\subsection{Phase field approach}

Consider for  $0 < \epsilon < \frac{2 \delta}{\pi}$ the function
\begin{displaymath}
\rho(x,t):= g \left( \frac{\phi(x,t)}{\epsilon} \right),
\end{displaymath}
where $g \in C^{1,1}(\mathbb{R})$ is given by
\begin{displaymath}
g(r)=
\left\{
\begin{array}{cl}
\cos^2(r), & | r | \leq \frac{\pi}{2}, \\[2mm]
0, & |r| > \frac{\pi}{2}.
\end{array}
\right.
\end{displaymath}
Note that $\mbox{supp}[\rho(\cdot,t)] = \overline{U_{\frac{\epsilon \pi}{2}}(t)} \subset U_{\delta}(t)$. Furthermore, we obtain from
the definition of $\nabla_{\phi}$ and (\ref{velext}) 
\begin{eqnarray}
\nabla_{\phi} \rho & = & \frac{1}{\epsilon} g'(\frac{\phi}{\epsilon}) \nabla_{\phi} \phi =0, \label{rho1} \\
\partial_t^{\bullet} \rho & = &  \frac{1}{\epsilon} g'( \frac{\phi}{\epsilon})
\bigl( \phi_t + (\vel,\nabla \phi) \bigr)=0. \label{rho2}
\end{eqnarray}

The phase field function $\rho$ allows us to approximate the integration over a surface $\Gamma(t)$  in terms of a volume integral over the diffuse interface. More precisely,
for fixed $t \in [0,T]$, the coarea formula implies for $\eta \in L^1(\Omega)$
\begin{displaymath}
\int_{\Omega} \eta \, \rho(\cdot,t) \, | \nabla \phi(\cdot,t) | \, dx  =  \int_{- \frac{\epsilon \pi}{2}}^{ \frac{\epsilon \pi}{2}} g \bigl( \frac{s}{\epsilon} \bigr)
\int_{\lbrace \phi(\cdot,t)=s \rbrace} \eta \, d \mathcal H^n ds    \approx  \frac{\epsilon \pi}{2} \int_{\lbrace \phi(\cdot,t)=0 \rbrace} \eta \, d \mathcal H^n
\end{displaymath}
for small $\epsilon >0$, so that we can view $\frac{2}{\epsilon \pi} \int_{\Omega} \eta \, \rho(\cdot,t) \, | \nabla \phi(\cdot,t) | \, dx$ as an approximation of $\int_{\Gamma(t)} \eta \, d \mathcal H^n$. 
This formula 
explains the appearance of the weight $\rho(\cdot,t) \, | \nabla \phi(\cdot,t) |$ in subsequent volume integrals.

In what follows we shall make use of the following
continuity properties of $\rho$.

\begin{lemma} \label{auxiliary} Let $s,t \in [0,T]$ with $| s-t | < \frac{\pi}{4 c_2} \epsilon$, $c_2$ as in (\ref{gradbound}).  Then 
$\mbox{supp}[ \rho(\cdot,s)] \subset U_{\frac{3 \epsilon \pi}{4}}(t)$ and 
\begin{eqnarray} 
| \rho(\cdot,t) - \rho(\cdot,s) | & \leq &  C \frac{|t-s|}{\epsilon} \sqrt{\rho(\cdot,t)} + C \frac{(t-s)^2}{\epsilon^2} \chi_{U_{\frac{3\epsilon \pi}{4}}(t)} \quad \mbox{ in } \Omega;
\label{rhodif} \\
| \rho_t(\cdot,t) - \rho_t(\cdot,s) | & \leq & C \frac{|t-s|}{\eps^2} \chi_{U_{\frac{3\epsilon \pi}{4}}(t)} \quad \mbox{ in } \Omega. \label{rhotdif}
\end{eqnarray}
\end{lemma}
{\bf Proof:} Let $s, t \in [0,T]$ with $| s-t| < \frac{\pi}{4 c_2} \epsilon$ and $x \in \mbox{supp}[ \rho(\cdot,s)]= \overline{U_{\frac{\epsilon \pi}{2}}(s)}$. Using the
mean value theorem and (\ref{gradbound}) we then have
\begin{displaymath}
| \phi(x, t)| \leq | \phi(x,s)| +  | \phi_t(x,\xi) | \, | t - s|   \leq \frac{\epsilon \pi}{2}+  c_2 | t-s|  < \frac{3 \epsilon \pi}{4},
\end{displaymath} 
i.e. $x \in U_{\frac{3 \epsilon \pi}{4}}(t)$. In order to prove (\ref{rhodif}) and (\ref{rhotdif}) we first observe that it is enough to verify the estimates for
$x \in U_{\frac{3\epsilon \pi}{4}}(t)$ in view of what we have just shown. There exists $\xi$ between $s$ and $t$ such that
\begin{equation} \label{rdif}
 | \rho(x,t) - \rho(x,s) |  = | \rho_t(x,\xi) | \, | t-s| = \frac{1}{\epsilon} | \phi_t(x,\xi) | \, \left| g'\left(\frac{\phi(x,\xi)}{\eps}\right) \right|  \, | t-s|   \leq  \frac{c_2 |t-s|}{\epsilon}  \left| g'\left(\frac{\phi(x,\xi)}{\eps}\right) \right|  
\end{equation}
by (\ref{gradbound}). Furthermore, since 
\begin{displaymath}
g'(r)=
\left\{
\begin{array}{cl}
-2 \sin(r) \cos(r), & | r | \leq \frac{\pi}{2}, \\[2mm]
0, & |r| > \frac{\pi}{2}
\end{array}
\right.
\end{displaymath}
we see immediately that
\begin{equation}  \label{sigmaprop}
\displaystyle
 | g'(r) |  \leq     2 \sqrt{g(r)}, \quad   | g'(r) - g'(\tilde{r})|   \leq   2 | r - \tilde{r}|,  \quad r, \tilde{r} \in \mathbb{R}.  
\end{equation}
As a result,
\begin{displaymath}
\left| g'\left(\frac{\phi(x,\xi)}{\eps}\right) \right|  \leq  \left| g'\left(\frac{\phi(x,t)}{\eps}\right) \right| + \frac{2}{\epsilon} | \phi(x,\xi) - \phi(x,t)| 
 \leq  2 \sqrt{\rho(x,t)} + \frac{2 c_2 |t-s| }{\epsilon}.  
\end{displaymath}
Inserting this bound into (\ref{rdif}) yields (\ref{rhodif}). Finally, using again (\ref{sigmaprop}) and (\ref{gradbound})  we obtain for  $x \in U_{\frac{3\epsilon \pi}{4}}(t)$  
\begin{eqnarray*}
| \rho_t(x,t) - \rho_t(x,s) | & \leq  &  \frac{1}{\epsilon} \left| g'\left(\frac{\phi(x,t)}{\epsilon} \right) - g'\left(\frac{\phi(x,s)}{\epsilon} \right)  \right| \, | \phi_t(x,t) | + 
\frac{1}{\eps} \left| g'\left(\frac{\phi(x,s)}{\epsilon} \right) \right| \, | \phi_t(x,t) - \phi_t(x,s)|   \\
& \leq &   \frac{C}{\epsilon^2} | \phi(x,t) - \phi(x,s) | + \frac{C}{\eps} | \phi_t(x,t) - \phi_t(x,s) | \leq \frac{C}{\eps^2} | t-s|. 
\end{eqnarray*}
\qed

\subsection{Discretization}

Suppose that $u$ is a smooth solution of (\ref{surfeq}). It is shown in Lemma \ref{extpde} of the Appendix that its extension $u^e$ satisfies the
strictly parabolic PDE
\be
\partial_t^{\bullet} u^e + u^e \, \nabla_{\phi} \cdot \vel - \frac{1}{| \nabla \phi | } \nabla \cdot \bigl( | \nabla \phi | \, \nabla  u^e \bigr)  = f^e + \phi \, R \quad \mbox{ in } 
\mathcal U_{\delta,T}
\label{eq:ue2}
\ee
with a smooth function $R$ depending on $u$ and $\phi$. In order to associate with (\ref{eq:ue2}) a suitable variational formulation 
we adapt an idea from \cite{ESSW11}, which uses an Eulerian transport identity. More precisely, we infer  
with the help of Lemma 3 in \cite{DE10}, (\ref{rho2}) and (\ref{eq:ue2}) that for every $\eta \in H^1(\Omega)$
\begin{eqnarray}
\lefteqn{ \frac{d}{dt}\io \ue\eta \,  \rho \,  \mgp = \io \bigl( \partial_t^{\bullet}( \ue \eta \rho ) + \ue \eta \rho \, \nabla_{\phi} \cdot \vel \bigr) \mgp } \nono \\
& = & \io \eta \bigl( \partial_t^{\bullet} u^e + u^e \, \nabphi \cdot \vel \bigr) \rho \,  \mgp + \io  u^e \partial_t^{\bullet} \eta \,  \rho \, \mgp  \nonumber  \\
&=&\io \eta \, \nabla \cdot \bigl( \mgp \, \nabla\ue \bigr) \rho + \io \eta \bigl( \fe+\phi R \bigr) \rho \, \mgp  + \io u^e \partial_t^{\bullet} \eta \, \rho \, \mgp  \label{weak} \\
&=&-\io (\nabla\ue, \nabla \eta)  \rho \,  \mgp+\io \fe\eta \,  \rho \,  \mgp+\io \ue (\vel,\nabla\eta) \rho \,  \mgp +\io    \phi \,   R \, \eta \, \rho \,  \mgp. \nonumber
\end{eqnarray}
Here, the last equality follows from integration by parts together with the fact that $(\nabla u^e, \nabla \rho)=  \frac{1}{\epsilon} g'\left( \frac{\phi}{\epsilon} \right)
(\nabla u^e,\nabla \phi)=0$ in view of (\ref{zetang}). \\
Let us first discretize with respect to time and denote by $0=t_0 < t_1< \ldots< t_M = T$ a partioning of $[0,T]$ with time steps $\tau_m:=t_m - t_{m-1}$ and $\tau:= \max_{m=1,\ldots,M} \tau_m$.
For a function $f=f(x,t)$ we shall write $f^m(x)=f(x,t_m)$. Integrating (\ref{weak}) with respect to $t \in (t_{m-1},t_m)$ we obtain for $\eta \in H^1(\Omega)$
\begin{eqnarray}
\lefteqn{ \hspace{-1cm} 
 \io  u^{e,m} \eta \rhop \mgpp-\io u^{e,m-1} \eta \rhon \mgpn+\tint\io  (\nabla u^e, \nabla \eta) \rho \,  \mgp } \label{weak1}  \\
&  & - \tint\io  u^e (\vel,\nabla \eta) \rho \,  \mgp = \tint\io \fe \,  \eta \,  \rho \,  \mgp + \tint\io \phi \,  R \,  \eta \,  \rho \, \mgp. \nonumber 
\end{eqnarray}
Neglecting the remainder term involving $R$ we now use the above relation in order discretize in space and hence to define our numerical method. 
In what follows we assume that $\Omega$ is polyhedral and consider a family  $(\mathcal T_h)_{0 < h \leq h_0}$ of triangulations of $\Omega$ 
with mesh size $h= \max_{T \in \mathcal T_h}h_T, \; h_T=\mbox{diam}(T)$. We assume that the family is regular in the sense that there exists 
$\sigma>0$ with
\begin{equation}  \label{regular}
r_T \geq \sigma h_T \qquad \forall T \in \mathcal T_h \quad  \forall 0 < h \leq h_0,
\end{equation}
where $r_T$ is the radius of the largest ball contained in $T$. Let us denote by 
$\mathcal N_h$ the set of vertices of the triangulation $\mathcal T_h$. 
 In order
to formulate our scheme
we require  a second phase field function with a slightly larger support, namely 
\begin{displaymath}
\tilde{\rho}(x,t)=g \bigl( \frac{\phi(x,t)}{2 \epsilon} \bigr), \quad 0<\epsilon < \frac{\delta}{\pi}.
\end{displaymath}
For $0 \leq m \leq M$ we then define
\begin{displaymath}
\mathcal T^m_h:= \lbrace T \in \mathcal T_h \, | \, \tilde{\rho}^m(x)>0 \mbox{ for some } x \in T \cap \mathcal N_h \rbrace \quad \mbox{ and } \quad
D^m_h:= \bigcup_{T \in \mathcal T^m_h} T
\end{displaymath}
as well as  the finite element space
\begin{displaymath}
V^m_h:= \lbrace v_h \in C^0(D^m_h) \, | \, v_{h|T} \mbox{ is a linear polynomial on each } T \in \mathcal T^m_h \rbrace.
\end{displaymath}
We denote by $I^m_h: C^0(D^m_h) \rightarrow V^m_h$ the standard Lagrange interpolation operator, i.e. $[I^m_h f](x)=f(x), x \in D^m_h \cap \mathcal N_h$. Note that
$D^m_h = \mbox{supp} I^m_h \tilde{\rho}^m$.

\begin{lemma} \label{dmh} Suppose that 
\begin{equation} \label{hecond}
\displaystyle
h \leq \frac{\cos^2(\frac{3 \pi}{8})}{2 c_1} \epsilon, \, \tau \leq \frac{\cos^2(\frac{3 \pi}{8})}{2 c_2} \epsilon.
\end{equation}
Then \\[2mm]
a) $U_{\frac{3 \epsilon \pi}{4}}(t) \subset D^m_h \subset U_{\frac{3 \epsilon \pi}{2}}(s)$ for all $s, t \in [\max(t_{m-1},0),\min(t_{m+1},T)], 0 \leq m \leq M$; \\[2mm]
b) $[I^m_h \tilde{\rho}^m](x) \geq \frac{1}{2} \cos^2(\frac{3 \pi}{8}), \, x \in U_{\frac{3 \epsilon \pi}{4}}(t_m), 0 \leq m \leq M$. 
\end{lemma}
{\bf Proof:} a) Let $x  \in D^m_h$, so that there exists $y \in \mathcal N_h$ such that  $|y-x| \leq h$ and
$\tilde{\rho}^m(y)>0$. Hence $| \phi^m(y) | < \epsilon \pi$ and  the mean value theorem together with (\ref{gradbound}) yields for $s \in [\max(t_{m-1},0),\min(t_{m+1},T)]$
\begin{eqnarray*}
| \phi(x,s) | & \leq &  | \phi(x,s) - \phi^m(x) | + | \phi^m(x) - \phi^m(y) | + | \phi^m(y) | \\
& < & | \phi_t(x,\xi) | \,  | s-t_m| +  | \nabla \phi^m(\eta) | \, | x-y | + \epsilon \pi  \\
& \leq &  c_2 \tau + c_1 h  + \epsilon \pi \leq \cos^2(\frac{3 \pi}{8}) \epsilon + \epsilon \pi \leq  \frac{3 \epsilon \pi}{2}
\end{eqnarray*}
in view of (\ref{hecond}). Hence, $x \in U_{\frac{3 \epsilon \pi}{2}}(s)$.  Next, let $x \in U_{\frac{3 \epsilon \pi}{4}}(t)$ for some $t \in [\max(t_{m-1},0),\min(t_{m+1},T)]$. Then $\tilde{\rho}(x,t) \geq \cos^2(\frac{3 \pi}{8})$
and we obtain similarly as above 
\begin{eqnarray*}
[I^m_h \tilde{\rho}^m](x) &  \geq &  \tilde{\rho}(x,t)  - | \tilde{\rho}(x,t) - \tilde{\rho}^m(x)| - | \tilde{\rho}^m(x) - [I^m_h \tilde{\rho}^m](x) | \\
& \geq &  \cos^2(\frac{3 \pi}{8}) -   | \tilde{\rho}_t(x,\xi) | \,  | t - t_m|   -   h \max_{y \in \overline{U_{\delta}(t_m)}} | \nabla \tilde{\rho}^m(y) | \\
& \geq & \cos^2(\frac{3 \pi}{8}) - c_2 \frac{\tau}{2 \epsilon} - c_1  \frac{h}{ 2 \epsilon}  \geq \frac{1}{2} \cos^2(\frac{3 \pi}{8}). 
\end{eqnarray*}
In particular, $[I^m_h \tilde{\rho}^m](x)>0$, so that $x \in D^m_h$. Using the above inequality for $t=t_m$ implies b). \qed \\[3mm]

Our finite element approximation of (\ref{surfeq}), (\ref{surfinit}) now reads:
For $m=1,2,..,M$  find $u^m_h \in V^m_h$ such that for all $v_h \in V^m_h$
\begin{eqnarray}
\lefteqn{ 
\int_{\Omega}  u^m_h \, v_h \, \rho^m \, | \nabla \phi^m | -
\int_{\Omega}  u^{m-1}_h \, v_h \,  \rho^{m-1} \, | \nabla \phi^{m-1} | + \tau_m \, \int_{\Omega}    (\nabla u^m_h, \nabla v_h) \, \rho^m \, | \nabla \phi^m | } \label{eq:fea} \\
 & &  - \tau_m \, \int_{\Omega}  u^m_h \, (\vel^m,\nabla v_h) \, \rho^m \, | \nabla \phi^m |   + \gamma  \tau_m^2 \, \int_{\Omega} I^m_h \tilde{\rho}^m
 ( \nabla u^m_h,\nabla v_h) 
 =  \tau_m \, \int_{\Omega}  f^{e,m} \, v_h \, \rho^m \, | \nabla \phi^m |. \nono
\end{eqnarray}
Here, $u^0_h \in V^0_h$ is defined as an $L^2$ projection of $u_0^e(x):=u_0(\tilde{p}(x,0)), x \in U_{\delta}(0)$, more precisely
\begin{equation} \label{defu0h}
\displaystyle \int_{D^0_h} u^0_h \, v_h = \int_{D^0_h} u_0^e \, v_h \qquad \forall v_h \in V^0_h.
\end{equation}
Furthermore, $f^{e,m}(x):=f(\tilde{p}(x,t_m),t_m), x \in U_{\delta}(t_m), 1 \leq m \leq M$. The parameter $\gamma>0$ will be chosen in
such a way as to ensure existence and stability for the scheme, see Lemma \ref{exis} and  Theorem \ref{discstab} below.

\begin{remark} \label{remarkscheme}
a) Lemma \ref{dmh} a) implies that  $\mbox{supp} \rho^m, \mbox{supp} \rho^{m-1} \subset   D^m_h = \mbox{supp} I^m_h \tilde{\rho}^m$,
so  that all integrals appearing in (\ref{eq:fea}) are taken only over $D^m_h$. In particular, if  $f \equiv 0$ we see
from the choice $v_h \equiv 1$ on $D^m_h$ that the scheme is mass conserving in the sense that
\begin{displaymath}
\io  u^m_h \, \rho^m \, | \nabla \phi^m | = \io  u^0_h \, \rho^0 \, | \nabla \phi^0|, \qquad m=1,\ldots,M.
\end{displaymath}
b) The term $\gamma  \tau_m^2 \, \int_{\Omega} I^m_h \tilde{\rho}^m ( \nabla u^m_h,\nabla v_h)$ introduces artificial diffusion into the scheme and  will play a crucial role in our 
analyis. A different form of stabilization is used in \cite{ESSW11}, Section 2.5. \\[2mm]
c) Unlike the schemes introduced in \cite{ESSW11} our method is not fully practical because we assume that the integrals are evaluated exactly.
In Section \ref{numerics} we shall follow \cite{ESSW11} in using numerical integration to obtain a fully practical scheme. A nice feature
of the resulting method is that the evolution of the hypersurfaces is tracked in a simple way via the evaluation of $\rho$. 
\end{remark}

In what follows we shall be concerned with the existence, stability and error bounds for (\ref{eq:fea}). The extension of our analysis to
the fully practical method mentioned above is currently out of reach and left for future research. However, the test calculations
in Section \ref{numerics} show that the parameter choices suggested by the analysis  work well 
also for the fully practical scheme.

\begin{lemma} \label{conn} There exists $0 < h_1 \leq h_0$ such that $D^m_h$ is connected for all $0< h \leq h_1$ and
$0 \leq m \leq M$.
\end{lemma}
{\bf Proof.}
To begin, we remark that there exists $0<h_1 \leq h_0$ and $\mu>0$ only depending on  $\sigma, c_0,c_1, c_2$ such that
for every $a \in \mathcal N_h \cap \overline{U_{\delta}(t)}$ there exists a neighbour $b \in \mathcal N_h$ with 
\begin{equation} \label{neighbour}
| \phi(a,t) - \phi(b,t) | \geq \mu h_T \qquad \mbox{ where } a,b, \in T
\end{equation}
for all $t \in [0,T], 0 < h \leq h_1$. 
Since $\Gamma(t_m)$ is connected  it is sufficient to show that for every $y \in D^m_h$ there exists $z \in \Gamma(t_m)$ and a path in $D^m_h$ connecting $y$ to $z$. Let us fix 
$y \in D^m_h$, say $y \in T$, where $\tilde{\rho}^m(x)>0$ for some $x \in T \cap \mathcal N_h$. We assume w.l.o.g. that $0 < \phi^m(x) < \epsilon \pi$. In view of (\ref{neighbour})
there exists a neighbour $x_1 \in \mathcal N_h$ of $x$ such that $\phi^m(x_1) \leq \phi^m(x) - \mu h_{\tilde{T}}$, where $x,x_1 \in \tilde{T}$. If $\phi^m(x_1) \leq 0$ then there
is $z \in [x,x_1]$ with $\phi^m(z)=0$. Hence, $z \in \Gamma(t_m)$ and the union of the segments $[y,x]$ and $[x,z]$ is a path in $D^m_h$ connecting $y$ to $z$. If $\phi^m(x_1)>0$,
then $\tilde{\rho}^m(x_1)>0$ so that  $[x,x_1] \subset D^m_h$ and we may repeat the above argument with $x$ replaced by $x_1$ and so on, until we reach $\Gamma(t_m)$ in a finite number of steps. \qed \\[3mm]

\begin{lemma} \label{exis} (Existence) Let $0 < h \leq h_1$. There exists $\tau_0>0$ such that the scheme (\ref{eq:fea}) has a unique solution $u^m_h \in V^m_h$ provided that
$0 < \tau \leq \tau_0$. 
\end{lemma}
{\bf Proof:} Since (\ref{eq:fea}) is equivalent to solving a linear system with a quadratic coefficient matrix, it is sufficient to prove that the following
problem only has the trivial solution: find $u_h \in V^m_h$ such that for all $v_h \in V^m_h$
\begin{eqnarray*} 
\lefteqn{\hspace{-3.5cm} \int_{\Omega} u_h \, v_h \,  \rho^m \, | \nabla \phi^m | + \tau_m  \int_{\Omega}    (\nabla u_h, \nabla v_h) \, \rho^m \,
 | \nabla \phi^m | - \tau_m \int_{\Omega}  u_h \, 
(\vel^m,\nabla v_h) \,  \rho^m \, | \nabla \phi^m | } \\
& & + \gamma   \tau_m^2 \, \int_{\Omega} I^m_h \tilde{\rho}^m ( \nabla u_h,  \nabla v_h) =0.
\end{eqnarray*}
Inserting $v_h=u_h$  we infer
\begin{eqnarray*}
\lefteqn{
\int_{\Omega} (u_h)^2  \rho^m \,  | \nabla \phi^m | + \tau_m \int_{\Omega}  | \nabla u_h |^2 \, \rho^m \, | \nabla \phi^m | + \gamma  \tau_m^2 \, \int_{\Omega}  I^m_h \tilde{\rho}^m  | \nabla u_h |^2 } \\
& = & \tau_m  \int_{\Omega} u_h \,  (\vel^m,\nabla u_h) \,  \rho^m \, | \nabla \phi^m | \leq  \tau_m \, \max_{x \in \overline{U_{\delta}(t_m)}} | \vel^m(x) | \, \io  | u_h | \, | \nabla u_h | \, \rho^m \,
| \nabla \phi^m |  \\
& \leq & \frac{1}{2} \int_{\Omega}  (u_h)^2  \, \rho^m \, | \nabla \phi^m | + \frac{1}{2} \tau \,  \bigl( \max_{x \in \overline{U_{\delta}(t_m)}} | \vel^m(x) | \bigr)^2 \tau_m  \int_{\Omega} | \nabla u_h |^2  \, \rho^m \, | \nabla \phi^m |.
\end{eqnarray*}
If we choose $\tau_0>0$ so small that $\frac{1}{2} \tau \,  \bigl( \max_{x \in \overline{U_{\delta}(t_m)}} | \vel^m(x) | \bigr)^2 \leq 1$  we deduce that
\begin{displaymath}
\int_{\Omega}  (u_h)^2 \, \rho^m \, | \nabla \phi^m | = \int_{\Omega} I^m_h  \tilde{\rho}^m  | \nabla u_h |^2=0,
\end{displaymath}
which implies that  $u_h \equiv 0$ on $\Gamma(t_m)$ and $\nabla u_h \equiv 0$ in $D^m_h$. According to Lemma \ref{conn},  
$D^m_h$ is connected, so that we conclude that $u_h \equiv 0$. \qed

\section{Stability bound}

The following lemma will be useful in estimating $L^2$–-integrals that are not weighted by $\rho$.

\begin{lemma} \label{aux}
There exists $C\geq 0$ such that for $t \in [0,T]$:
\be
\int_{U_{\frac{3 \epsilon \pi}{4}}(t)}  f^2 \leq C \io  f^2 \rho(\cdot,t) | \nabla \phi(\cdot,t)|   +C\eps^2\int_{U_{\frac{3 \epsilon \pi}{4}}(t)}  | \nabla f |^2 \qquad \mbox{ for all } f \in H^1(\Omega).
\ee
\end{lemma}
\begin{remark}
Note that Lemma \ref{dmh} b) implies that
\begin{equation}  \label{stab1}
\displaystyle \int_{U_{\frac{3 \epsilon \pi}{4}}(t_m)}  | \nabla f |^2 \leq \frac{2}{\cos^2(\frac{3 \pi}{8})} \io I^m_h \tilde{\rho}^m | \nabla f |^2, \quad f \in H^1(\Omega), m=0,\ldots,M.
\end{equation}
\end{remark}
{\bf Proof.} We may assume that $f$ is smooth, the general case then follows with the help of an approximation argument. Since $F_t$ is a diffeomorphism from $\Gamma(0) \times
(-\frac{3 \epsilon \pi}{4}, \frac{3 \epsilon \pi}{4})$ onto $U_{\frac{3 \epsilon \pi}{4}}(t)$, the transformation rule yields
\begin{equation} \label{aux1}
\displaystyle 
c_1 \int_{U_{\frac{3\epsilon \pi}{4}}(t)} f(x)^2 dx  \leq  \int_{-\frac{3 \eps\pi}{4}}^{\frac{3 \eps\pi}{4}}\int_{\Gamma(0)} f(F_t(P,s))^2 do_P ds \leq c_2 \int_{U_{\frac{3\epsilon \pi}{4}}(t)} f(x)^2 dx.
\end{equation}
The definition of $F_t$ together with (\ref{odesys}) implies for $|s|\leq \frac{3 \eps\pi}{4}$, $|\tilde{s} | \leq \frac{\eps\pi}{4}$
\been
f(F_t(P,s))&=& f(F_t(P,\tilde{s}))+ \int_{\tilde{s}}^s \bigl( \nabla f(F_t(P,r)), \frac{\pd F_t}{\pd r}(P,r) \bigr) dr \\
&=& f(F_t(P,\tilde{s}))+\int_{\tilde{s}}^s \bigl( \nabla f(F_t(P,r)), \frac{\nabla \phi(F_t(P,r),t)}{|\nabla \phi(F_t(P,r),t)|^2} \bigr) dr
\eeen
and therefore
\begin{eqnarray*}
f(F_t(P,s))^2 & \leq &  2  f(F_t(P,\tilde{s}))^2 + C \eps \int_{-\frac{3 \eps\pi}{4}}^{\frac{3 \eps\pi}{4}} | \nabla f (F_t(p,r)) |^2 dr  \\
& \leq & C  f(F_t(P,\tilde{s}))^2 \rho(F_t(P,\tilde{s}),t)   +C\eps \int_{-\frac{3 \eps\pi}{4}}^{\frac{3 \eps\pi}{4}}  | \nabla f (F_t(p,r)) |^2 dr,
\end{eqnarray*}
since $\rho(F_t(P,\tilde{s}),t) = \cos^2 \bigl(\frac{\phi(F_t(P,\tilde{s}),t)}{\epsilon} \bigr) = \cos^2(\frac{\tilde{s}}{\epsilon}) \geq \cos^2(\frac{\pi}{4}), | \tilde{s} | \leq
\frac{\epsilon \pi}{4}$. 
Integrating with respect to $P \in \Gamma(0), s \in (-\frac{3 \eps\pi}{4},\frac{3 \eps\pi}{4})$ and recalling (\ref{aux1}) we obtain for 
$|\tilde{s} | \leq \frac{\eps\pi}{4}$
\begin{eqnarray*}
 \int_{U_{\frac{3 \epsilon \pi}{4}}(t)}  f(x)^2 dx & \leq &     C \epsilon \int_{\Gamma(0)} f(F_t(P,\tilde{s}))^2 \rho(F_t(P,\tilde{s}),t) do_P  + C \epsilon^2 
\int_{-\frac{3 \eps\pi}{4}}^{\frac{3 \eps\pi}{4}}  \int_{\Gamma(0)}  |  \nabla f (F_t(p,r)) |^2 do_P dr   \\  
& \leq & C \epsilon \int_{\Gamma(0)}    f(F_t(P,\tilde{s}))^2 \rho(F_t(P,\tilde{s}),t) do_P  + C\eps^2 \int_{U_{\frac{3 \epsilon \pi}{4}}(t)}  | \nabla f(x) |^2 dx.
\end{eqnarray*}
If we integrate with respect to $\tilde{s} \in (-\frac{\eps\pi}{4},\frac{\eps\pi}{4})$, divide by $\epsilon$ and recall (\ref{gradbound}) we obtain the assertion. \qed  

It follows from Theorem 4.4 in \cite{DE07} (extended in a straightforward way to the case of a nontrivial $f$) that
(\ref{surfeq}), (\ref{surfinit}) has a unique solution $u$ which satisfies
\begin{displaymath}
\sup_{(0,T)} \Vert u(\cdot,t) \Vert_{L^2(\Gamma(t))}^2 + \int_0^T \Vert \nabla_{\Gamma} u(\cdot,t) \Vert_{L^2(\Gamma(t))}^2 dt
\leq c \bigl( \Vert u_0 \Vert_{L^2(\Gamma(0))}^2 + \int_0^T \Vert f(\cdot,t) \Vert_{L^2(\Gamma(t))}^2 dt \bigr).
\end{displaymath}
The following theorem gives a discrete version of this estimate in the phase field setting.

\begin{theorem} \label{discstab}
Suppose that (\ref{hecond}) holds.
There exist $\gamma_1 >0$ and  $\tau_1 \leq \tau_0$ such that  
\begin{eqnarray*}
\lefteqn{ \hspace{-2cm}
\max_{m=1,\ldots,M} \frac{2}{\epsilon \pi} \int_{\Omega}   (u^m_h)^2 \, \rho^m \, | \nabla \phi^m | + \, \sum_{m=1}^M  \tau_m \, 
\frac{2}{\epsilon \pi} \int_{\Omega}  | \nabla u^m_h |^2   \rho^m \, | \nabla \phi^m | }\\
& \leq &    C \bigl(  \int_{\Gamma(0)} (u_0)^2 +  \sum_{m=1}^M \tau_m \int_{\Gamma(t_m)} (f^m)^2 \bigr),
\end{eqnarray*}
provided that $\gamma \geq \gamma_1$   and $\tau \leq \max \bigl( \tau_1, \epsilon^2 \bigr) $.
\end{theorem}
{\bf Proof:} Setting $v_h=\up$ in (\ref{eq:fea}) 
we find after a straighforward calculation
\begin{eqnarray}
\lefteqn{ \hspace{-1cm}
 \frac12\io (\up)^2 \, \rhop \, \mgpp - \frac12\io(\un)^2 \, \rhon \, \mgpn+\frac12\io (\up-\un)^2 \, \rhon \, \mgpn }\nono\\
& & \quad + \tau_m  \io |\nabla\up|^2 \, \rhop \, \mgpp + \gamma  \tau_m^2 \io I^m_h \tilde{\rho}^m  | \nabla \up |^2 \nono \\
& = &  - \frac12\io (\up)^2 \, (\rhop-\rhon) \, \mgpn + \frac12\io (\up)^2 \, \rhop \, (\mgpn-\mgpp)\nono\\
& & \quad + \tau_m \io \up (\vp,\nabla\up) \, \rhop \,  \mgpp+  \tau_m \io f^{e,m} \up \, \rhop \, \mgpp \nono \\
&:= &  I + II +  III + IV.  \label{eq:aa}
\end{eqnarray}
Clearly,
\be
I  = - \frac12\int_{t_{m-1}}^{t_m} \io   (\up)^2 \, \rho_t(\cdot,s) \, \mgpn ds,
\label{eq:a}
\ee
while 
\begin{eqnarray} 
II & = &- \frac12\tau_m \io (\up)^2 \, \rhop \,  (\nabla \phi_t^m, \nu^m)  + \frac{1}{2} \io (\up)^2 \, \rhop \, \bigl( 
\mgpn - \mgpp +\tau_m (\nabla \phi^m_t,\nu^m) \bigr) \nono  \\
 & = &  II_1 + II_2.
\label{eq:b}
\end{eqnarray}
Integrating by parts and abbreviating
$H^m = - \nabla \cdot \nu^m$ we obtain
\begin{eqnarray}
II_1 &=& \frac12 \tau_m \io (\up)^2 \, (\nabla\rhop, \nu^m) \, \phi_t^m +  \tau_m \io \up (\nabla\up,\nu^m) \rhop \,  \phi_t^m 
+ \frac12 \tau_m \io (\up)^2 \, \nabla \cdot \nu^m \, \rhop \,  \phi_t^m \nono\\
&=& \frac12 \tau_m \io (\up)^2 \, \rho_t^m \, \mgpp +  \tau_m \io \up (\nabla\up,\nu^m) \rhop \, \phi_t^m  \nono 
 - \frac12 \tau_m  \io (\up)^2 H^m \, \rhop \,  \phi_t^m, \label{eq:tt}
\end{eqnarray}
since
\begin{displaymath}
( \nabla \rho^m, \nu^m) \phi^m_t = \frac{1}{\epsilon} g'(\frac{\phi^m}{\epsilon}) \, \phi^m_t \, ( \nabla \phi^m,\nu^m)  
= \rho^m_t \, | \nabla \phi^m |.
\end{displaymath}
In order to rewrite $III$ we first observe that in view of (\ref{tanggrad}) and  (\ref{velext})
\begin{displaymath}
(\vp, \nabla \up)  =  (\vp, \nabla_{\phi^m} \up)  + (\nabla \up, \nu^m) (\vp, \nu^m)  =  (\vp, \nabla_{\phi^m} \up) - (\nabla \up, \nu^m) \frac{ \phi^m_t} {\mgpp},
\end{displaymath}
so that (\ref{Eulintparts}),  (\ref{rho1}) and again (\ref{velext}) imply
\begin{eqnarray}
III & = & \frac12 \tau_m \io (\vp, \nabla_{\phi^m}(\up)^2) \rhop \, \mgpp- \tau_m \io \up (\nabla\up, \nu^m) \rhop \,  \phi_t^m   \nono\\
&=& - \frac12 \tau_m  \io  \nabla_{\phi^m } \cdot \vp (\up)^2 \, \rhop \, \mgpp 
 - \frac12 \tau_m \io H^m ( \vp, \nu^m)  (\up)^2 \, \rhop \, \mgpp \nono \\
 & & - \tau_m  \io \up (\nabla\up, \nu^m) \rhop \,  \phi_t^m \label{eq:bb}\\
&=& - \frac12 \tau_m  \io  \nabla_{\phi^m} \cdot \vp  (\up)^2 \, \rhop \, \mgpp + \frac12 \tau_m \io (\up)^2H^m \, \rhop \,  \phi_t^m  - \tau_m \io \up (\nabla\up, \nu^m) \, \rhop \,  \phi_t^m . \nono
\end{eqnarray}
Inserting  (\ref{eq:a})-(\ref{eq:bb}) into (\ref{eq:aa}) we infer that 
\begin{eqnarray}
\lefteqn{ \hspace{-1cm}
 \frac12\io (\up)^2 \, \rhop \,\mgpp - \frac12\io (\un)^2 \, \rhon \, \mgpn  
 + \tau_m \io |\nabla\up|^2 \, \rhop \, \mgpp + \gamma  \tau_m^2 \io I^m_h \tilde{\rho}^m | \nabla \up |^2 \nono }  \\
& \leq  & \frac12\int_{t_{m-1}}^{t_m}\io (\up)^2  \left(\rho_t^m \mgpp-\rho_t(\ldot,s)\mgpn\right)  -\frac12 \tau_m \io \nabla_{\phi^m}\cdot\vp (\up)^2 \, \rhop \, \mgpp \nono\\
& & + \frac{1}{2} \io (\up)^2 \, \rhop \, \bigl( 
\mgpn - \mgpp +\tau_m (\nabla \phi^m_t,\nu^m) \bigr)
+ \tau_m  \io  f^{e,m} \up \, \rhop \,\mgpp.  \label{ap} 
\end{eqnarray}
We deduce from (\ref{rhotdif}),  Lemma \ref{aux}, (\ref{stab1}) and the assumption $\tau \leq \epsilon^2$ that
\begin{eqnarray*}
\lefteqn{ \hspace{-5mm} 
 \left| \frac{1}{2}  \tint\io (\up)^2 (\rho^m_t\mgpp-\rho_t(\ldot,s)\mgpn) \right|  } \\
 & \leq & C \tint \io   (\up)^2 \bigl(  | \rho^m_t - \rho_t(\ldot,s) | + | \rho_t(\ldot,s) |  \, | \mgpp -\mgpn | \bigr) \\
 &\leq & 
 C \frac{\tau_m^2}{\epsilon^2} \int_{U_{\frac{3 \epsilon \pi}{4}}(t_m)} (\up)^2   
   \leq  C  \frac{\tau_m^2}{\epsilon^2} \io  (\up)^2 \, \rhop \mgpp  + C \tau_m^2 \int_{U_{\frac{3 \epsilon \pi}{4}}(t_m)} | \nabla \up |^2 \\
 & \leq &   C \tau_m \io  (\up)^2 \,  \rhop \mgpp  + (\gamma-1)  \tau_m^2 \io I_h \tilde{\rho}^m \, | \nabla \up |^2
\end{eqnarray*}
if we choose $\gamma \geq \gamma_1:=C+1$. Finally, using Taylor expansion and (\ref{gradbound}) we infer that
\begin{eqnarray*}
\lefteqn{ \hspace{-1.5cm}
 \left|  -\frac12 \tau_m  \io \nabla_{\phi^m}\cdot\vp(\up)^2 \, \rhop \, \mgpp  + \frac{1}{2} \io (\up)^2 \, \rhop \, \bigl( 
\mgpn - \mgpp +\tau_m (\nabla \phi^m_t,\nu^m) \bigr) \right|  } \\
& \leq &   C \tau_m   \io  (\up)^2 \, \rhop \mgpp + C \tau_m^2 \int (\up)^2 \rho^m \leq C \tau_m \io (\up)^2 \rho^m \mgpp.
\end{eqnarray*}
Inserting the above estimates into (\ref{ap}) we find
\begin{eqnarray}
\lefteqn{ \hspace{-5mm}
 \frac12\io (\up)^2 \, \rhop \, \mgpp  + \tau_m  \io |\nabla\up|^2 \, \rhop \, \mgpp +  \tau_m^2 \io I^m_h \tilde{\rho}^m | \nabla \up |^2 }  \label{stab2}  \\
& \leq & \frac12\io (\un)^2 \, \rhon \, \mgpn +  C \tau_m \io  (\up)^2 \, \rho^m \, \mgpp  + \tau_m \io  (f^{e,m})^2 \, \rho^m \, \mgpp. \nono
\end{eqnarray}
If $\tau_1 \leq \tau_0$ is sufficiently small we therefore deduce for $\tau \leq \tau_1$
\begin{eqnarray*}
\lefteqn{ \hspace{-2cm}
\io (\up)^2 \, \rhop \, \mgpp + \tau_m  \io |\nabla\up|^2 \, \rhop \, \mgpp } \\
& \leq &  (1+ C \tau_m)  \io (\un)^2 \, \rhon \, \mgpn + C \tau_m  \io  (f^{e,m})^2 \, \rho^m \,  \mgpp,
\end{eqnarray*}
from which we obtain after summation from $m=1,\ldots,l$ and division by $\epsilon$ that
\begin{eqnarray*}
\lefteqn{ \hspace{-5mm} 
\frac{1}{\epsilon} \io  (u_h^l)^2 \, \rho^l \,  | \nabla \phi^l | + \sum_{m=1}^l \tau_m \frac{1}{\epsilon} \io |\nabla\up|^2 \, \rhop \, \mgpp } \\
& \leq & \frac{1}{\epsilon} \io  (u_h^0)^2  \, \rho^0 \, |  \nabla \phi^0 |  + C \sum_{m=0}^{l-1} \tau_{m+1} \frac{1}{\epsilon} \io  (\up)^2 \, \rhop \, \mgpp +
C \sum_{m=1}^l \tau_m \frac{1}{\epsilon}  \io  (f^{e,m})^2 \, \rho^m \, \mgpp.
\end{eqnarray*}
Using Lemma \ref{dmh} a), (\ref{defu0h}) and (\ref{zeest1}) we may estimate
\begin{displaymath}
\frac{1}{\epsilon} \io  (u_h^0)^2 \, \rho^0 \,  | \nabla \phi^0 | \leq \frac{C}{\epsilon} \int_{D^0_h} (u^0_h)^2 \leq \frac{C}{\epsilon} \int_{D^0_h} (u^e_0)^2 \leq \frac{C}{\epsilon} 
\int_{U_{\frac{3 \epsilon \pi}{2}}(0)} (u^e_0)^2
\leq C \int_{\Gamma(0)} (u_0)^2.
\end{displaymath}
Arguing in a similar way for the term involving $f^{e,m}$ we derive
\begin{eqnarray}
\lefteqn{  \hspace{-5mm}
\frac{1}{\epsilon} \io  (u_h^l)^2 \, \rho^l \, | \nabla \phi^l | + \sum_{m=1}^l \tau_m \frac{1}{\epsilon} \io |\nabla\up|^2 \, \rhop \, \mgpp } \nono \\
& \leq &   C \sum_{m=0}^{l-1} \tau_{m+1} \frac{1}{\epsilon} \io  (\up)^2 \, \rhop \, \mgpp 
+ C \bigl( \int_{\Gamma(0)} (u_0)^2  +  \sum_{m=1}^l \tau_m \int_{\Gamma(t_m)} (f^m)^2 \bigr). \label{eq:dd}
\end{eqnarray} 
The discrete Gronwall inequality yields the bound on 
$\max_{m=1,\ldots,M} \frac{1}{\epsilon} \int_{\Omega}   (u^m_h)^2 \, \rho^m \, | \nabla \phi^m |$, which combined with (\ref{eq:dd}) implies the
second inequality. \qed

\section{Error estimate}

Before we formulate our error bound we derive interpolation estimates that are adapted to our setting.

\begin{lemma} \label{interpol} Suppose that (\ref{hecond}) holds and let $z^e$ be defined by (\ref{ext}). Then we have for $m=1,\ldots,M$ and $t \in [t_{m-1},t_m]$:
\begin{eqnarray*}
\int_{D^m_h} | (z^e - I^m_h z^e)(\cdot,t) |^2 + h^2 \int_{D^m_h} | \nabla (z^e - I^m_h z^e)(\cdot,t) |^2 \leq C \epsilon h^4 \Vert z(\cdot,t) \Vert_{H^2(\Gamma(t))}^2, \\
\hspace{-3cm}  \int_{D^m_h} | (z^e_t - I^m_h z^e_t)(\cdot,t) |^2 \leq C \epsilon h^4 \bigl( \Vert  \partial_t^{\bullet} z(\cdot,t) \Vert_{H^2(\Gamma(t))}^2 + 
\Vert z(\cdot,t) \Vert_{H^3(\Gamma(t))}^2 \bigr).
\end{eqnarray*}
\end{lemma}
{\bf Proof:} Let $t \in [t_{m-1},t_m]$. Standard interpolation theory together with Lemma \ref{dmh} a) and (\ref{zeest1}) implies that
\begin{eqnarray*}
\lefteqn{
\int_{D^m_h} | (z^e - I^m_h z^e)(\cdot,t) |^2  + h^2 \int_{D^m_h} | \nabla (z^e - I^m_h z^e)(\cdot,t) |^2 }  \\
& \leq & c h^4  \int_{D^m_h}  | D^2 z^e(\cdot,t) |^2 \leq c h^4
\int_{U_{\frac{3 \epsilon \pi}{2}}(t)} | D^2 z^e(\cdot,t) |^2   
 \leq  C \epsilon  h^4 \Vert z(\cdot,t) \Vert_{H^2(\Gamma(t))}^2.
\end{eqnarray*}
The second bound follows in the same way using (\ref{zeest2}). \qed \\[3mm]

\begin{theorem} \label{thm:main}Suppose that the solution of (\ref{surfeq}), (\ref{surfinit}) satisfies
\begin{equation} \label{ureg}
\displaystyle
\max_{t \in [0,T]} \Vert u(\cdot,t) \Vert_{H^2(\Gamma(t))}^2 + \int_0^T \bigl( \Vert u(\cdot,t) \Vert_{H^3(\Gamma(t))}^2 + \Vert \partial_t^{\bullet} u(\cdot,t) \Vert_{H^2(\Gamma(t))}^2 \bigr) dt < \infty.
\end{equation}
Then there exists $0< \tau_2 \leq \tau_1$ and a constant $C\geq 0$ such that 
\begin{displaymath}
\max_{m=1,\ldots,M} \frac{2}{\epsilon \pi} \int_{\Omega}  \,  | u^{e,m} - u^m_h|^2 \, \rho^m \,  | \nabla \phi^m | +  \sum_{m=1}^M  \tau_m
\frac{2}{\epsilon \pi} \int_{\Omega}  | \nabla (u^{e,m}-   u^m_h) |^2 \rho^m \, | \nabla \phi^m | \leq   C \epsilon^2,
\end{displaymath}
provided that $\tau \leq \max(\epsilon^2,\tau_2)$, $\gamma \geq \gamma_1$ and (\ref{hecond}) hold. 
\end{theorem}

{\bf Proof:} 
Let us write
\begin{displaymath}
u^{e,m} - u^m_h = (u^{e,m} - I^m_h u^{e,m}) + (I^m_h u^{e,m} - u^m_h)=: d^m + e^m_h.
\end{displaymath}
If we combine (\ref{weak1}) for $\eta = v_h \in V^m_h$ with (\ref{eq:fea})  we find
\begin{eqnarray*}
\lefteqn{ 
\io \ep v_h \,  \rhop  \mgpp-\io \en v_h \,   \rhon \mgpn+ \tau_m \io (\nabla\ep,\nabla v_h) \rhop \mgpp }\\
& & -\tau_m  \io \ep (\vp,\nabla v_h) \rhop  \mgpp+ \gamma \tau_m^2\io I^m_h \tilde{\rho}^m (\nabla\ep,\nabla v_h) \\
& = & \left[ - \io  \dhp v_h \, \rhop  \mgpp + \io \dhn v_h \, \rhon  \mgpn \right]  -  \tau_m  \io  (\nabla\dhp,\nabla v_h) \rhop \mgpp \\
& & + \tau_m  \io \dhp (\vp,\nabla v_h) \rhop  \mgpp + \gamma \tau_m^2\io I^m_h \tilde{\rho}^m (\nabla I_h u^{e,m},\nabla v_h) \\
& & + \tint\io\left[  (\nabla u^{e,m}, \nabla v_h) \rhop  \mgpp - (\nabla u^e,\nabla v_h)  \rho \,  \mgp \right]\\
& & + \tint\io\left[ u^e (\vel, \nabla v_h) \rho \, \mgp-  u^{e,m} (\vp, \nabla v_h) \rhop  \mgpp\right] \\
&  & +\tint\io\left[\fe v_h \, \rho \,  \mgp-  f^{e,m} v_h \, \rhop \mgpp\right] +\tint\io \phi \,  R \,  v_h \, \rho \, \mgp \\
& =: &  \sum_{i=1}^8 \langle S_i^m, v_h \rangle.
\end{eqnarray*}
Inserting $v_h= e^m_h$ and following the argument in the proof of Theorem \ref{discstab} leading to (\ref{stab2}) we obtain
\begin{eqnarray}
\lefteqn{ 
 \frac12\io (\ep)^2 \rhop \mgpp + \tau_m  \io  |\nabla\ep|^2 \rhop \mgpp +   \tau_m^2 \io I^m_h \tilde{\rho}^m | \nabla \ep | ^2 } \nono \\
& \leq & \frac12\io (\en)^2 \rhon \mgpn +  C \tau_m \io  (\ep)^2 \rho^m \mgpp  + \sum_{i=1}^8 \langle S^m_i,\ep \rangle. \label{eq:cc} 
\end{eqnarray}
We now deal individually with the terms $\langle S^m_i,\ep \rangle, i =1,\ldots,8$ in (\ref{eq:cc}). 
Clearly,
\begin{eqnarray*}
\lefteqn{ \hspace{-2cm} 
| \langle  S_1^m, \ep \rangle | \leq C \io  | d^m - d^{m-1} | \, | \ep |  \,  \rhop  + C \io    | d^{m-1} | \,  | \ep | \,  | \nabla ( \phi^m - \phi^{m-1}) | 
\, \rhop } \\
& & \qquad + C \io   | d^{m-1} | \, | \ep | \, | \rhop - \rhon |  \equiv I + II + III.
\end{eqnarray*}
In order to estimate $I$ we first deduce from Lemma \ref{dmh} a) that every $T \in \mathcal T_h$ with $T \cap \mbox{supp} \rho^m \neq \emptyset$ satisfies $T \in \mathcal T^{m-1}_h \cap \mathcal T^m_h$.
Therefore $I^{m-1}_h u^{e,m-1} = I^m_h u^{e,m-1}$ on $\mbox{supp} \rho^m$, which yields
\begin{displaymath}
d^m - d^{m-1} = [u^{e,m}- u^{e,m-1}]- I^m_h[u^{e,m}- u^{e,m-1}] = \int_{t_{m-1}}^{t_m} (u^e_t - I^m_h u^e_t)(\cdot,t)  \quad \mbox{ on } \mbox{supp} \rho^m. 
\end{displaymath}
Hence, Lemma \ref{dmh} a), Lemma \ref{interpol} and (\ref{gradbound}) imply that
\begin{eqnarray*}
| I | & \leq & C \io \int_{t_{m-1}}^{t_m} | u^e_t - I^m_h u^e_t | \,  | \ep | \,  \rhop \leq C \sqrt{\tau_m}  \left( \io  (\ep)^2 \rhop \right)^{\frac{1}{2}} \left( \int_{t_{m-1}}^{t_m}   
\int_{D^m_h}  | u^e_t - I^m_h u^e_t |^2 \right)^{\frac{1}{2}} \\
& \leq & \tau_m \io (e^m_h)^2 \rhop \mgpp + C \epsilon h^4 \int_{t_{m-1}}^{t_m} \bigl( \Vert  \partial^{\bullet}_t u(\cdot,t) \Vert_{H^2(\Gamma(t))}^2 + 
\Vert u(\cdot,t) \Vert_{H^3(\Gamma(t))}^2 \bigr) dt
\end{eqnarray*}
and similarly,
\begin{eqnarray*}
| II | & \leq &   C \tau_m \io  | d^{m-1} | \, | \ep | \, \rhop \leq C \tau_m \left( \io (\ep)^2  \rhop \right)^{\frac{1}{2}} \left( \int_{D^{m-1}_h} | d^{m-1} |^2 
\right)^{\frac{1}{2}} \\
& \leq & \tau_m \io  (e^m_h)^2 \rhop \mgpp + C \epsilon h^4 \tau_m  \Vert u^{m-1} \Vert_{H^2(\Gamma(t_{m-1}))}^2.
\end{eqnarray*}
Next, we deduce from (\ref{rhodif}), Lemma \ref{dmh} a),  (\ref{gradbound}),  Lemma \ref{interpol}, Lemma \ref{aux}  and (\ref{stab1}) that
\begin{eqnarray*}
\lefteqn{
| III |  \leq   C \frac{\tau_m}{\epsilon} \io  | d^{m-1} \, | \ep | \, \sqrt{\rhop} + C \frac{\tau_m^2}{\epsilon^2} \int_{U_{\frac{3 \epsilon \pi}{4}}(t_m)} | d^{m-1} | \, | \ep | } \\
& \leq & \tau_m \io (\ep)^2 \rho^m + C \frac{\tau_m}{\epsilon^2} \Vert d^{m-1} \Vert_{L^2(D^{m-1}_h)}^2 + C \frac{\tau_m^2}{\epsilon^2} \left( \int_{ U_{\frac{3 \epsilon \pi}{4}}(t_m)} ( \ep )^2  
\right)^{\frac{1}{2}} \Vert d^{m-1} \Vert_{L^2(D^{m-1}_h)} \\
& \leq & C \tau_m \io (\ep)^2 \rho^m \mgpp  + C \frac{\tau_m h^4}{\epsilon}  \Vert u^{m-1} \Vert_{H^2(\Gamma(t_{m-1}))}^2  \\
& & \quad + C \frac{\tau_m^2 h^2}{\epsilon^{\frac{3}{2}}} \Vert u^{m-1} \Vert_{H^2(\Gamma(t_{m-1}))}  \left(
 \io  ( \ep )^2 \rhop \mgpp + \epsilon^2  \int_{ U_{\frac{3 \epsilon \pi}{4}}(t_m)}
| \nabla e^m_h |^2 \right)^{\frac{1}{2}}   \\
& \leq & C \tau_m \io (e^m_h)^2 \rhop \mgpp + \frac{\tau_m^2}{8} \io I^m_h  \tilde{\rho}^m | \nabla e^m_h |^2 + C \frac{\tau_m h^4}{\epsilon} \Vert u^{m-1} \Vert_{H^2(\Gamma(t_{m-1}))}^2,
\end{eqnarray*}
where we used that $\tau \leq \epsilon^2$. Again by Lemma \ref{interpol} we have
\begin{eqnarray*} 
| \langle S^m_2,\ep \rangle | & \leq &  \tau_m \left( \io  | \nabla \ep |^2 \rhop \mgpp \right)^{\frac{1}{2}} \left( \int_{D^m_h} | \nabla d^m |^2 \right)^{\frac{1}{2}} \\
& \leq &  \frac{1}{8} \tau_m \io   | \nabla \ep |^2 \rhop \mgpp  + C \tau_m \epsilon h^2 \Vert u^m \Vert_{H^2(\Gamma(t_m))}^2,
\end{eqnarray*}
while
\begin{eqnarray*}
 | \langle S^m_3,\ep \rangle | & \leq &  C \tau_m \io  | d^m | \, | \nabla \ep | \, \rhop \mgpp  \leq  C \tau_m  \left( \io  | \nabla \ep |^2 \rhop \mgpp \right)^{\frac{1}{2}} 
\left( \int_{D^m_h} | d^m|^2 \right)^{\frac{1}{2}}   \\
& \leq &  \frac{1}{8} \tau_m \io   | \nabla \ep |^2 \rhop \mgpp + C \tau_m \epsilon  h^4  \Vert u^m \Vert_{H^2(\Gamma(t_m))}^2.
\end{eqnarray*}
Lemma \ref{dmh} a), (\ref{zeest1}) and Lemma \ref{interpol} yield
\begin{eqnarray*}
 | \langle S^m_4,\ep \rangle | & \leq &  C \tau_m^2 \left( \io I^m_h \tilde{\rho}^m |  \nabla \ep |^2 \right)^{\frac{1}{2}} \left( \int_{D^m_h} | \nabla I^m_h u^{e,m} |^2 \right)^{\frac{1}{2}} \\
 & \leq & \frac{\tau_m^2}{8} \io I^m_h \tilde{\rho}^m | \nabla \ep |^2 + C \tau_m^2 \int_{D^m_h} \bigl( | \nabla u^{e,m} |^2 + | \nabla d^m |^2 \bigr) \\
 & \leq & \frac{\tau_m^2}{8} \io I^m_h \tilde{\rho}^m | \nabla \ep |^2 + C \tau_m^2 \epsilon \Vert u^m \Vert_{H^2(\Gamma(t_m))}^2.
\end{eqnarray*}
We deduce from (\ref{rhodif}), Lemma \ref{dmh} a) and Lemma \ref{extest} that
\begin{eqnarray*}
\lefteqn{ 
| \langle S^m_5,\ep \rangle |  \leq   C \int_{t_{m-1}}^{t_m} \io \left[ | \nabla (u^{e,m} - u^e) | \, \rho^m  + | \nabla u^e | \, | \nabla (\phi^m - \phi) | \, \rho^m 
 +  | \nabla u^e | \, | \rho^m - \rho | \right] | \nabla \ep| } \\
& \leq & C \tau_m  \int_{t_{m-1}}^{t_m} \io  | \nabla u^e_t | \, | \nabla \ep | \rho^m + C \tau_m  \int_{t_{m-1}}^{t_m} \io  | \nabla u^e | \, | \nabla \ep | \rho^m  \\
& &  \quad + C \frac{\tau_m}{\epsilon}  \int_{t_{m-1}}^{t_m} \io | \nabla u^e | \, | \nabla \ep | \sqrt{\rho^m} + C \frac{\tau_m^2}{\epsilon^2}  \int_{t_{m-1}}^{t_m}  \int_{U_{\frac{3 \epsilon \pi}{4}}(t_m)}
| \nabla u^e | \, | \nabla \ep| \\
& \leq & C \left[  \tau_m^{\frac{3}{2}} \left( \int_{t_{m-1}}^{t_m} \Vert u^e_t(\cdot,t) \Vert^2_{H^1(U_{\frac{3 \epsilon \pi}{4}}(t))} \right)^{\frac{1}{2}} 
+  \frac{\tau_m^2}{\epsilon} \max_{t_{m-1} \leq t \leq t_m} \Vert u^e(\cdot,t) \Vert_{H^1(U_{\frac{3 \epsilon \pi}{4}}(t))} \right] \left( \io | \nabla \ep|^2 \rho^m \right)^{\frac{1}{2}} \\
& & \quad  + C \frac{\tau_m^3}{\epsilon^2} \max_{t_{m-1} \leq t \leq t_m} \Vert u^e(\cdot,t) \Vert_{H^1(U_{\frac{3 \epsilon \pi}{2}}(t))} \left( \int_{U_{\frac{3 \epsilon \pi}{4}}(t_m)} | \nabla \ep |^2 \right)^{\frac{1}{2}} \\
& \leq & \frac{\tau_m}{8} \io | \nabla \ep |^2 \rho^m \mgpp + C \tau_m^2 \epsilon \int_{t_{m-1}}^{t_m} \bigl( \Vert \partial_t^{\bullet} u(\cdot,t) \Vert_{H^1(\Gamma(t))}^2 + \Vert u(\cdot,t) \Vert_{H^2(\Gamma(t))}^2
\bigr) dt \\
& & + \frac{\tau_m^2}{8} \io I^m_h \tilde{\rho}^m | \nabla \ep |^2 + C \tau_m^2 \epsilon \max_{t_{m-1} \leq t \leq t_m} \Vert u(\cdot,t) \Vert_{H^1(\Gamma(t))}^2.
\end{eqnarray*}
Here we have used again that $\tau_m \leq \tau \leq \epsilon^2$. In a similar way we obtain
\begin{eqnarray*}
\lefteqn{ \hspace{-1cm} 
| \langle S^m_6,\ep \rangle |  \leq    C \int_{t_{m-1}}^{t_m} \io  \left[ |  u^{e,m} - u^e| \, \rho^m  + | u^e | \, | \vel^m | \nabla \phi^m |  - \vel | \nabla \phi |  | \, \rho^m 
 +  |  u^e | \, | \rho^m - \rho | \right] | \nabla \ep| } \\
& \leq & \frac{\tau_m}{8} \io | \nabla \ep |^2 \rho^m \mgpp + C \tau_m^2 \epsilon \int_{t_{m-1}}^{t_m} \bigl( \Vert \partial_t^{\bullet} u(\cdot,t) \Vert_{L^2(\Gamma(t))}^2 + \Vert u(\cdot,t) \Vert_{H^1(\Gamma(t))}^2
\bigr) dt \\
& & + \frac{\tau_m^2}{8} \io I^m_h \tilde{\rho}^m | \nabla \ep |^2 + C \tau_m^2 \epsilon \max_{t_{m-1} \leq t \leq t_m} \Vert u(\cdot,t) \Vert_{L^2(\Gamma(t))}^2
\end{eqnarray*}
as well as
\begin{eqnarray*}
| \langle S^m_7,\ep \rangle | & \leq & C \int_{t_{m-1}}^{t_m} \io \left[  | f^e | \nabla \phi | - f^{e,m} | \nabla \phi^m | \, | \,  | \ep | \, \rhop + | f^{e,m} | \,  | \ep |  \, | \rho - \rhop | \right]  \\
& \leq &  C \tau_m^2 \io  | \ep | \, \rhop + C \frac{\tau_m^2}{\epsilon} \io  | \ep | \, \sqrt{\rhop} + C \frac{\tau_m^3}{\epsilon^2} \int_{U_{\frac{3 \epsilon \pi}{4}}(t_m)} | \ep |  \\
& \leq & \tau_m \io  ( \ep )^2 \rho^m + C \frac{\tau_m^3}{\epsilon}   + C \frac{\tau_m^3}{\epsilon^{\frac{3}{2}}}  \left( \io ( \ep )^2 \rho^m + \epsilon^2 \int_{U_{\frac{3 \epsilon \pi}{4}}(t_m)}
| \nabla \ep |^2 \right)^{\frac{1}{2}} \\
& \leq & C \tau_m \io  (\ep)^2 \rhop \mgpp + \frac{\tau_m^2}{8}  \io I^m_h \tilde{\rho}^m | \nabla \ep |^2 + C \frac{\tau_m^3}{\epsilon},
\end{eqnarray*}
where we have used that $| U_{\frac{3 \epsilon \pi}{4}}(t_m) | \leq C \epsilon$  and again the fact that $\tau \leq \epsilon^2$. Finally, since $R$ is bounded and $| \phi(\cdot,t) | \leq c \epsilon$ on $\mbox{supp} \rho(\cdot,t), \, t \in [t_{m-1},t_m]$ we may
estimate with the help of  (\ref{rhodif}) and Lemma \ref{aux}
\begin{eqnarray*}
| \langle S^m_8,\ep \rangle | & \leq & C \int_{t_{m-1}}^{t_m} \io \left[  | \phi | \, | \ep | \, \rhop +  | \phi | \, | \ep | \, | \rho - \rhop |  \right]  \\
& \leq & C \epsilon \tau_m \io  | \ep | \, \rhop + C \tau_m^2 \io  | \ep | \, \sqrt{\rhop} + C \frac{\tau_m^3}{\epsilon} \int_{U_{\frac{3 \epsilon \pi}{4}}(t_m)} | \ep | \\
& \leq & \tau_m \io ( \ep )^2 \rho^m  + C \tau_m \epsilon^3 + C \tau_m^3 \epsilon + C \frac{\tau_m^3}{\sqrt{\epsilon}} \left( \io ( \ep )^2 \rho^m + \epsilon^2 \int_{U_{\frac{3 \epsilon \pi}{4}}(t_m)}
| \nabla \ep |^2 \right)^{\frac{1}{2}} \\
& \leq & C \tau_m \io  (\ep)^2 \rhop \mgpp  + \frac{\tau_m^2}{8} \io I^m_h \tilde{\rho}^m | \nabla \ep |^2 + C \tau_m \epsilon^3 + C \tau_m^3 \epsilon.
\end{eqnarray*}
Inserting the above estimates into (\ref{eq:cc}) we obtain 
\begin{eqnarray*}
\lefteqn{ 
 \frac12\io (\ep)^2 \rhop \mgpp + \frac{\tau_m}{2} \io |\nabla\ep|^2 \rhop \mgpp + \frac{\tau_m^2}{4} \io I^m_h \tilde{\rho}^m | \nabla \ep |^2 } 
    \\
& \leq & \frac12\io (\en)^2 \rhon \mgpn +  C \tau_m \io  (\ep)^2 \rho^m \mgpp + C \bigl( \frac{\tau_m^3}{\epsilon}  + \tau_m \epsilon^3 \bigr) \nonumber  \\
& & +  C \tau_m \, \epsilon \bigl( h^2 + \frac{h^4}{\epsilon^2} + \tau \bigr) \max_{t_{m-1} \leq t \leq t_m} \Vert u(\cdot,t) \Vert_{H^2(\Gamma(t))}^2 \nonumber  \\
&& + C \epsilon \bigl( h^4 + \tau^2 \bigr) \int_{t_{m-1}}^{t_m} \bigl( \Vert  \partial^{\bullet}_t u(\cdot,t) \Vert_{H^2(\Gamma(t))}^2 + 
\Vert u(\cdot,t) \Vert_{H^3(\Gamma(t))}^2 \bigr)dt. \nono
\end{eqnarray*}
Choosing $\tau_2 \leq \tau_1$ small enough and using (\ref{hecond}) as well as  $\tau \leq \epsilon^2$  we infer
\begin{eqnarray*}
\lefteqn{ \hspace{-2cm} 
\io (\ep)^2 \rhop \mgpp + \tau_m \io |\nabla\ep|^2 \rhop \mgpp } \\
& \leq &   (1+ C \tau_m)  \io (\en)^2 \rhon \mgpn 
 + C \epsilon^3 \tau_m \max_{t_{m-1} \leq t \leq t_m} \Vert u(\cdot,t) \Vert_{H^2(\Gamma(t))}^2  \\
 & & + C \epsilon^5 \int_{t_{m-1}}^{t_m} \bigl( \Vert  \partial^{\bullet}_t u(\cdot,t) \Vert_{H^2(\Gamma(t))}^2 + 
\Vert u(\cdot,t) \Vert_{H^3(\Gamma(t))}^2 \bigr) dt + C  \tau_m \epsilon^3. 
\end{eqnarray*}
Summing from $m=1,\ldots,l$, dividing by $\epsilon$ and recalling (\ref{ureg}) we derive
\begin{eqnarray*}
\lefteqn{ \hspace{ -1cm} 
\ove\io (e^l_h)^2 \rho^l |\nabla\phi^l| + \sum_{m=1}^l \tau_m \frac{1}{\epsilon} \io |\nabla\ep|^2 \rhop \mgpp } \\
& \leq & \ove \io (e^0_h)^2\rho^0 \, | \nabla \phi^0 | +
 \sum_{m=0}^{l-1} \tau_{m+1} \ove \io  (\ep)^2 \rhop \,  \mgpp + C \epsilon^2.
\end{eqnarray*}
In order to estimate the first term on the right hand side we write $e^0_h = (I^0_h u^e_0 - u^e_0) + (u^e_0- u_h^0)$ and recall the definition (\ref{defu0h}) of $u^0_h$ as an $L^2$ projection:
\begin{displaymath}
\io  (e^0_h)^2 \rho^0 \, | \nabla \phi^0 | \leq C \int_{D^0_h} (e^0_h)^2 \leq C \int_{D^0_h} | u^e_0 - I^0_h u^e_0 |^2 \leq C \epsilon h^4 \Vert u_0 \Vert_{H^2(\Gamma(0))}^2
\end{displaymath}
by Lemma \ref{interpol}. Thus 
\begin{equation} \label{err}
\displaystyle 
\ove\io  (e^l_h)^2 \rho^l |\nabla\phi^l| + \sum_{m=1}^l \tau_m \frac{1}{\epsilon} \io |\nabla\ep|^2 \rhop \mgpp \leq  \sum_{m=0}^{l-1} \tau_{m+1} \ove \io  (\ep)^2 \rhop \mgpp + C \epsilon^2 
\end{equation}
and  the discrete Gronwall lemma gives
\begin{equation}  \label{err1}
\displaystyle 
\max_{m=1,\ldots,M} \ove\io (\ep)^2 \rhop \mgpp \leq C \epsilon^2.
\end{equation}
The remainder of the proof  follows from (\ref{err}) and  Lemma \ref{interpol}. \qed

Using the result of Theorem \ref{thm:main} we can now also derive an error bound on the surface.

\begin{corollary} \label{cormain} In addition to the assumptions of Theorem \ref{thm:main} suppose that $\displaystyle \max_{t \in [0,T]} \Vert u(\cdot,t) \Vert_{W^{2,\infty}(\Gamma(t))} < \infty$. 
Furthermore we assume that there exists $\alpha>0$ such that $h_T \geq \alpha \epsilon$ for all $T \in \mathcal T_h$ with $| T \cap \Gamma(t) |>0, t \in [0,T]$. Then
\begin{displaymath}
\max_{m=1,\ldots,M} \int_{\Gamma(t_m)} | u^m - u^m_h |^2 + \sum_{m=1}^M  \tau_m \int_{\Gamma(t_m)} | \nabla_{\Gamma} ( u^m - u^m_h ) |^2 \leq C \epsilon^2.
\end{displaymath}
\end{corollary}
{\bf Proof:} Let us fix $m \in \lbrace 1,\ldots,M \rbrace$ and define $\mathcal T^m_{\Gamma,h}:= \lbrace T \in \mathcal T_h \, | \, | T \cap \Gamma(t_m) |>0 \rbrace$. Hence, given $T \in \mathcal T^m_{\Gamma,h}$, there
exists $x_T \in \Gamma(t_m)$ with $\phi^m(x_T)=0$. We infer from (\ref{gradbound}) and (\ref{hecond}) that for arbitrary $x \in T$ 
\begin{displaymath}
| \phi^m(x) | = | \phi^m(x) - \phi^m(x_T) |  \leq c_1 | x - x_T| \leq c_1 h_T \leq \frac{\epsilon}{2} \cos^2 \left( \frac{3 \pi}{8} \right) \leq \frac{\epsilon \pi}{4},
\end{displaymath}
and therefore
\begin{equation}  \label{gammalow}
\rho^m(x) \geq \frac{1}{2} \quad \mbox{ for all } x \in T, \, T \in \mathcal T^m_{\Gamma,h}.
\end{equation}
We now argue in a similar way as in \cite{DDEH10}, page 368. Using an interpolation inequality and an inverse estimate we infer that
\begin{eqnarray*}
\int_{\Gamma(t_m)} | u^m - u^m_h |^2 & = & \sum_{T \in \mathcal T^m_{\Gamma,h}} \int_{T \cap \Gamma(t_m)} | u^m - u^m_h |^2 \leq 2  \sum_{T \in \mathcal T^m_{\Gamma,h}} | T \cap \Gamma(t_m) | \left(
\Vert d^m \Vert_{L^{\infty}(T)}^2 + \Vert e^m_h \Vert_{L^{\infty}(T)}^2 \right)  \\
& \leq & C \sum_{T \in \mathcal T^m_{\Gamma,h} } | T \cap \Gamma(t_m) | \, h_T^2 \Vert \nabla u^{e,m} \Vert_{W^{1,\infty}(T)}^2 + C \sum_{T \in \mathcal T^m_{\Gamma,h} } h_T^n  h_T^{-(n+1)} \Vert e^m_h \Vert_{L^2(T)}^2  \\
& \leq & C h^2 | \Gamma(t_m) | \Vert u^m \Vert_{W^{1,\infty}(\Gamma(t_m))}^2 + C \epsilon^{-1} \sum_{T \in \mathcal T^m_{\Gamma,h} } \int_T | e^m_h |^2 \rho^m \, | \nabla \phi^m |,
\end{eqnarray*}
where the last inequality follows from (\ref{gammalow}), (\ref{gradbound}) and the assumption that $h_T \geq \alpha \epsilon, T \in \mathcal T^m_{\Gamma,h}$.
In a similar way we obtain
\begin{displaymath}
\int_{\Gamma(t_m)} | \nabla_{\Gamma} (u^m-u^m_h) |^2  \leq C h^2 | \Gamma(t_m) | \Vert u^m \Vert_{W^{2,\infty}(\Gamma(t_m))}^2  + C \epsilon^{-1} \sum_{T \in \mathcal T^m_{\Gamma,h} } 
\int_T |  \nabla e^m_h |^2 \rho^m \, | \nabla \phi^m |.
\end{displaymath}
Thus,
\begin{eqnarray*}
\lefteqn{ \hspace{-1.5cm}  \max_{m=1,\ldots,M} \int_{\Gamma(t_m)} | u^m - u^m_h |^2 + \sum_{m=1}^M  \tau_m \int_{\Gamma(t_m)} | \nabla_{\Gamma} ( u^m - u^m_h ) |^2  \leq 
C h^2  \max_{t \in [0,T]} \Vert u(\cdot,t) \Vert_{W^{2,\infty}(\Gamma(t))}^2 } \\
& & \quad \quad + C \epsilon^{-1} \max_{m=1,\ldots,M} \int_{\Omega} | e^m_h |^2 \rho^m | \nabla \phi^m |
+ C \epsilon^{-1} \sum_{m=1}^M \tau_m \int_{\Omega} | \nabla e^m_h |^2 \rho^m | \nabla \phi^m | \\
& \leq & C \epsilon^2,
\end{eqnarray*}
by (\ref{hecond}), (\ref{err1}) and (\ref{err}). \qed

\section{Numerical Results} \label{numerics}

As already mentioned in Remark \ref{remarkscheme} c), the scheme (\ref{eq:fea}), (\ref{defu0h}) is not fully practical. Therefore, our implementation uses the
following modification: Find $u^m_h \in V^m_h, m=0,1,\ldots,M$ such that
\begin{eqnarray}
\lefteqn{ \hspace{-1cm}
\int_{\Omega}  u^m_h \, v_h \, I^m_h \rho^m \, | \nabla I^m_h \phi^m | -
\int_{\Omega}  u^{m-1}_h \, v_h \, I^{m-1}_h \rho^{m-1} \, | \nabla I^{m-1}_h  \phi^{m-1} | } \label{eq:feam} \\
&& + \tau_m \, \int_{\Omega}    (\nabla u^m_h, \nabla v_h) \, I^m_h \rho^m \, | \nabla I^m_h  \phi^m | 
- \tau_m \, \int_{\Omega}  u^m_h \, (I_h^m \hat{\vel}^m,\nabla v_h) \,  I^m_h \rho^m \, | \nabla I^m_h \phi^m | \nono \\
& & + \gamma  \tau_m^2 \, \int_{\Omega} I^m_h \tilde{\rho}^m ( \nabla u^m_h,\nabla v_h)  
 =  \tau_m \, \int_{\Omega}  I_h^m\hat{f}^m \, v_h \, I^m_h  \rho^m \, | \nabla I^m_h  \phi^m | \nono
\end{eqnarray}
for all  $v_h \in V^m_h$ and $1 \leq m \leq M$. Here, $\hat{\vel}^m(x):= \vel(\hat{p}(x,t_m),t_m)$, $\hat{f}^m(x)=
f(\hat{p}(x,t_m),t_m)$, where $\hat{p}(x,t)$ denotes the closest point projection of a point $x$ onto $\Gamma(t)$.
Setting $\hat{u}_0(x)=u_0(\hat{p}(x,0))$ we define the initial data $\hat{u}^0_h \in V^0_h$ by
\begin{equation} \label{defu0hm}
\displaystyle \int_{D^0_h} \hat{u}^0_h \, v_h = \int_{D^0_h} I^0_h \hat{u}_0 \, v_h \qquad \forall v_h \in V^0_h.
\end{equation}
Let us remark that the evaluation of  
$\hat{p}(x,t)$ is  easier compared to $\tilde{p}(x,t)$, which has been used to extend the data for the scheme (\ref{eq:fea}), (\ref{defu0h}). However, we claim that
\begin{equation} \label{pdif}
\tilde{p}(x,t) - \hat{p}(x,t) = O(\phi(x,t)^2).
\end{equation}
To see this, we first observe that $\hat{p}(x,t)$ is characterized by the conditions
\begin{displaymath}
\phi(\hat{p}(x,t),t)=0 \quad \mbox{ and } \quad x - \hat{p}(x,t) \perp \Gamma(t) \mbox{ at } \hat{p}(x,t).
\end{displaymath}
Therefore, it is not difficult to verify with the help of Taylor expansion that
\begin{displaymath}
x - \hat{p}(x,t) = \lambda(x,t) \nabla \phi(\hat{p}(x,t),t), \quad \mbox{ with } \lambda(x,t) = \frac{\phi(x,t)}{| \nabla \phi(\hat{p}(x,t),t)|^2} + O(\phi(x,t)^2).
\end{displaymath}
Combining this relation with (\ref{pk}) in the Appendix we find that
\begin{displaymath}
\tilde{p}(x,t) - \hat{p}(x,t) = \phi(x,t) \left[ \frac{\nabla \phi(\hat{p}(x,t),t)}{| \nabla \phi(\hat{p}(x,t),t) |^2} - \frac{ \nabla \phi(x,t)}{ | \nabla \phi(x,t) |^2} \right]
+ O(\phi(x,t)^2) = O(\phi(x,t)^2).
\end{displaymath}
In particular, we infer from (\ref{pdif}) that replacing $\tilde{p}$ by $\hat{p}$ in the extension of $\vel,f$ and $u_0$ will not affect the result of 
Theorem \ref{thm:main}. In contrast, it is not straightforward to handle the interpolation terms $I^m_h \rho^m$ and $I^{m-1}_h \rho^{m-1}$ in (\ref{eq:feam}).
Applying a standard interpolation estimate to $\rho^m - I^m_h \rho^m$ will result in a term of the form $h^2 \Vert \rho^m \Vert_{H^2} \approx \frac{h^2}{\epsilon^2}$, which
we are currently not able to analyze. 
The results of our test calculations below however show that the use of the interpolation operator in  (\ref{eq:feam}), (\ref{defu0hm}) does not lead to reduced convergence rates. 
More
precisely we investigate the experimental order of convergence (eoc) for the following errors:
\begin{eqnarray*}
\mathcal{E}_1 &= &\max_{m=1,\ldots,M} \frac{2}{\epsilon \pi} \int_{\Omega}  \,  | I_h^m\hat{u}^m - u^m_h|^2 \, I^m_h \rho^m\, | \nabla I^m_h \phi^m |, \\
\mathcal{E}_2 & = & \frac{2}{\epsilon \pi} \sum_{m=1}^M  \tau_m \int_{\Omega}  | \nabla (I_h^m \hat{u}^m-   u^m_h) |^2 I^m_h \rho^m \, | \nabla I^m_h \phi^m |,
\end{eqnarray*}
where $\hat{u}^m(x)=u(\hat{p}(x,t_m),t_m)$. We use the finite element toolbox Alberta 2.0, \cite{alberta}, and implement a similar mesh refinement strategy to that in \cite{BNS} with a fine mesh constructed in $D_h^m$ 
and a coarser mesh in $\Omega \backslash D_h^m$. The linear systems appearing in each time step were solved using GMRES together with diagonal preconditioning.
The values of $h$ given below are such that $h:=\max_{T\in D_h^m}h_T$, $h_T=\mbox{diam}(T)$.


%
%
%
%
%


\subsection{2D examples}\label{2d}

We set $\Omega = (-2.4, 2.4)^2$, $T=0.1$,  and choose $\gamma=0.01, \, \epsilon=85.33 \, h$ as well as a uniform time step $\tau_m = 0.0025\eps^2, m=1,\ldots,M$.
In all our examples below $\Gamma(t)$ will be a circle $\Gamma(t)= \lbrace x \in \mathbb{R}^2 \, | \, | x - m(t) |= 1 \rbrace$ of radius 1 with
center $m(t) \in \mathbb{R}^2$. In addition to $\mathcal E_1, \mathcal E_2$ we shall also investigate 
the errors appearing in Corollary \ref{cormain}. To do so we choose $L>0$ and define the following quadrature points
\begin{displaymath}
x_l(t):= m(t)+  \bigl(\cos(\frac{2 \pi l}{L}), \sin(\frac{2 \pi l}{L}) \bigr)^T, \quad l=0,\ldots,L-1
\end{displaymath}
as well as
\begin{eqnarray*}
\mathcal{E}_3 & = & \max_{m=1,\ldots,M} \sum_{l=0}^{L-1}\frac{2\pi}{L}| u(x_l(t_m),t_m) - u^{m}_h(x_l(t_m)) |^2,\\
\mathcal{E}_4 & = & \sum_{m=1}^M  \tau_m\sum_{l=0}^{L-1}\frac{2\pi}{L}| \nabla_{\Gamma}u(x_l(t_m),t_m) - \nabla_{\Gamma} u^{m}_h(x_l(t_m)) |^2.
\end{eqnarray*}
In our computations $L=200$ turned out to be sufficient. \\

{\it Example 1} The form of our first example is similar to Example 7.2 in \cite{DE07}. We
consider $\Gamma(t)=\Gamma=S^1, t \in [0,T]$ described as the zero level set of the function $\phi(x):= x_1^2+x_2^2-1$. The function $u:S_T \rightarrow \mathbb{R}, u(x,t):= \frac{1}{2} e^{-4t}(x_1^2-x_2^2)$
is a solution of (\ref{surfeq}), (\ref{surfinit}) with $\vel=\mathbf{0}, f=0$ and initial data $u_0(x) = \frac{1}{2}(x_1^2-x_2^2)$.   
Due to the symmetry of the problem we only solve on $\Omega = (0, 2.4)^2$ and impose homogenous Neumann boundary conditions on the symmetry boundaries. 
In Table \ref{tab:ex1} we display the values of $\mathcal{E}_i$, $i=1\to4$, together with the eocs. We see that the eoc for $\mathcal{E}_1$ is reducing towards $4$, the eocs for for $\mathcal{E}_2$ and $\mathcal{E}_3$ are close to $4$ and the eoc for $\mathcal{E}_4$ is between 2 and 3.


 \begin{center}
\begin{table}[!h]
\resizebox{\textwidth}{!}{
 \begin{tabular}{ |c|c||c|c|c|c|c|c|c|c| }
 \hline
$h$ & $\varepsilon$ & $\mathcal{E}_1$ & $eoc_1$ & $\mathcal{E}_2$ & $eoc_2$ & $\mathcal{E}_3$ & $eoc_3$ & $\mathcal{E}_4$ & $eoc_4$  \\ 
 \hline
 \hline
{4.6875e-03} & $0.4$ & {3.8512e-05} & - & {2.3461e-04} & - & {1.7447e-05} & - & {3.0051e-06} & -   \\ 
{3.3146e-03} & $0.2\sqrt{2}$ & {6.1934e-06} & {5.273} & {5.8552e-05} & {4.005}   & {5.3838e-06} & {3.393}  & {1.3270e-06} & {2.358}  \\ 
{2.3437e-03} & $0.2$ & {1.1872e-06} & {4.766} & {1.4660e-05} & {3.996}  & {1.4812e-06}  & {3.724}  & {5.5173e-07} & {2.532}  \\ 
{1.6573e-03} & $0.1\sqrt{2}$ & {2.5553e-07} & {4.432} & {3.6713e-06} & {3.995}  & {3.8766e-07} & {3.868}  & {2.4519e-07} & {2.340}  \\ 
{1.1719e-03} & $0.1$ & {5.8950e-08} & {4.232} & {9.2008e-07} & {3.993}  & {9.9099e-08}  & {3.936}  & {9.8614e-08} & {2.628}  \\ 
\hline
\end{tabular}}
\caption{Errors and experimental orders of convergence for Example 1 \label{tab:ex1}}
\end{table}
\end{center}

{\it Example 2}
We again consider the stationary unit  circle $\Gamma(t)=\Gamma = S^1$ together with the same level set function as in the previous example. 
The function $u(x,t):= e^{-4t} \left[ x_1 x_2 \cos(\pi t) + \frac{1}{2}(x_1^2 - x_2^2) \sin(\pi t) \right]$ 
is a solution of  (\ref{surfeq}), (\ref{surfinit}) for the velocity field
$\bbv ({x}) = \frac{\pi}{2} (x_2, -x_1)^T, f=0$ and the initial data $u_0(x)=x_1 x_2$. A similar choice of velocity appears in Example 3 in
\cite{DE10}.
The results are displayed in Table \ref{tab:ex2} where we see eocs that are very similar to the ones in Table \ref{tab:ex1}. 

 \begin{center}
\begin{table}[!h]
\resizebox{\textwidth}{!}{
 \begin{tabular}{ |c|c||c|c|c|c|c|c|c|c| }
 \hline
$h$ & $\varepsilon$ & $\mathcal{E}_1$ & $eoc_1$ & $\mathcal{E}_2$ & $eoc_2$ & $\mathcal{E}_3$ & $eoc_3$ & $\mathcal{E}_4$ & $eoc_4$  \\ 
 \hline
 \hline
{4.6875e-03} & $0.4$ & {2.0565e-04} & - & {1.0763e-03} & - & {2.7651e-05} & - & {4.3137e-06} & -   \\ 
{3.3146e-03} & $0.2\sqrt{2}$ & {3.2822e-05} & {5.295} & {2.7030e-04} & {3.987}   & {8.1077e-06} & {3.540}  & {1.6031e-06} & {2.856}  \\ 
{2.3437e-03} & $0.2$ & {6.5608e-06} & {4.645} & {6.7864e-05} & {3.988}  & {2.1848e-06}  & {3.784}  & {5.9541e-07} & {2.858}  \\ 
{1.6573e-03} & $0.1\sqrt{2}$ & {1.4513e-06} & {4.353} & {1.7017e-05} & {3.991}  & {5.6637e-07} & {3.895}  & {2.3962e-07} & {2.626}  \\ 
{1.1719e-03} & $0.1$ & {3.4022e-07} & {4.186} & {4.2668e-06} & {3.991}  & {1.4412e-07}  & {3.949}  & {9.6590e-08} & {2.622}  \\ 
\hline
\end{tabular}}
\caption{Errors and experimental orders of convergence for Example 2  \label{tab:ex2}}
\end{table}
\end{center}

{\it Example 3} (cf. \cite[Section 3.1]{ESSW11}, \cite{TLLWV09}, Example 5.2)
We consider the family of unit circles $\Gamma(t) = \lbrace x \in \mathbb{R}^2 \, | \, (x_1 + \frac{1}{2} -2t)^2 + x_2^2 =1 \rbrace$
described as the zero level set of $\phi(x,t)=(x_1 + \frac{1}{2}-2t)^2+x_2^2-1$. The function $u:S_T \rightarrow \mathbb{R}$,
$u(x,t)= e^{-4t}(x_1 + \frac{1}{2}-2t) x_2$ is a solution of (\ref{surfeq}), (\ref{surfinit}) for the 
velocity field $\vel (x,t) = (2,0)^T,f=0$ and the inital data $u_0(x)=(x_1+\frac{1}{2})x_2$.
The results are displayed in Table \ref{tab:ex4} where we see eocs that are similar to the ones in Tables \ref{tab:ex1} and \ref{tab:ex2}. 

 \begin{center}
\begin{table}[!h]
\resizebox{\textwidth}{!}{
 \begin{tabular}{ |c|c||c|c|c|c|c|c|c|c| }
 \hline
$h$ & $\varepsilon$ & $\mathcal{E}_1$ & $eoc_1$ & $\mathcal{E}_2$ & $eoc_2$ & $\mathcal{E}_3$ & $eoc_3$ & $\mathcal{E}_4$ & $eoc_4$  \\ 
 \hline
 \hline
{4.6875e-03} & $0.4$ & {1.5537e-04} & - & {9.3201e-04} & - & {1.8431e-05} & - & {3.0082e-06} & -   \\ 
{3.3146e-03} & $0.2\sqrt{2}$ & {2.5206e-05} & {5.248} & {2.3280e-04} & {4.002}   & {5.6312e-06} & {3.421}  & {1.2489e-06} & {2.537}  \\ 
{2.3437e-03} & $0.2$ & {4.8726e-06} & {4.742} & {5.8500e-05} & {3.985}  & {1.5443e-06}  & {3.733}  & {4.8015e-07} & {2.758}  \\ 
{1.6573e-03} & $0.1\sqrt{2}$ & {1.0558e-06} & {4.413} & {1.4776e-05} & {3.970}  & {4.0396e-07} & {3.869}  & {1.9389e-07} & {2.616}  \\ 
{1.1719e-03} & $0.1$ & {2.4507e-07} & {4.214} & {3.7865e-06} & {3.929}  & {1.0350e-07}  & {3.929}  & {8.1747e-08} & {2.492}  \\ 
\hline
\end{tabular}}
\caption{Errors and experimental orders of convergence for Example 3  \label{tab:ex4}}
\end{table}
\end{center}

\subsection{3D example}
{\it Example 4}
Here we consider the first example in Section 7 of \cite{LNR11} in which a family of expanding  and collapsing spheres is considered such that 
$\Gamma(t)= \lbrace x \in \mathbb{R}^2 \, | \, |x|=r(t) \rbrace$ where $r(t) = 1 + \sin^2(\pi t)$, 
described as the zero level set of $\phi(x,t)=x_1^2+x_2^2+x_3^2-r(t)^2$. The function $u:S_T \rightarrow \mathbb{R}, u(x,t)= \frac{2}{r(t)^2|x|^2} e^{-6\int_0^t\frac1{r^2(t)}} x_1 x_3$
is a solution of (\ref{surfeq}), (\ref{surfinit}) for the velocity field $\bbv(x,t) = \frac{r'(t)}{|x|} x, f=0$ and 
the initial data $u_0(x)=\frac{2}{|x|^2}x_1 x_3$. 
We set $\Omega=(-4,4)^3$ and choose $\gamma=0.01, \, \epsilon = 1.85 \, h$  as well as  a uniform time step $\tau_m = 0.5h^2, m=1,\ldots,M$.
For this example we only display the errors on the surfaces which are in this case approximated by
the quadrature rules
\begin{eqnarray*}
\mathcal{E}_3 & = & \max_{m=1,\ldots,M} \sum_{k=0}^{2L-1} \sum_{l=0}^{L-1} (\frac{\pi}{L})^2 | u(x_{k,l}(t_m),t_m) - u^{m}_h(x_{k,l}(t_m)) |^2 \, \sin(\frac{l \pi}{L}),\\
\mathcal{E}_4 & = & \sum_{m=1}^M  \tau_m \sum_{k=0}^{2L-1} \sum_{l=0}^{L-1} ( \frac{\pi}{L} )^2 | \nabla_{\Gamma}u(x_{k,l}(t_m),t_m) - \nabla_{\Gamma} u^{m}_h(x_{k,l}(t_m)) |^2 \, \sin(\frac{l \pi}{L}),
\end{eqnarray*}
where
\begin{displaymath}
x_{k,l}(t)= r(t) \bigl( \cos(\frac{k \pi}{L}) \sin( \frac{l \pi}{L}), \sin ( \frac{k \pi}{L}) \sin( \frac{l \pi}{L}), \cos( \frac{l \pi}{L}) \bigr)^T, \quad k=0,\ldots,2L-1, l=0,\ldots,L-1.
\end{displaymath}
For the choice $L=200$ the results are displayed in Table \ref{tab:ex3d}, where we see eocs close to $4$ for $\mathcal{E}_3$ and eocs close to $2$ for $\mathcal{E}_4$. 


 \begin{center}
\begin{table}[!h]
 \begin{tabular}{ |c|c||c|c|c|c| }
 \hline
$h$ & $\varepsilon$ & $\mathcal{E}_3$ & $eoc_3$ & $\mathcal{E}_4$ & $eoc_4$  \\
 \hline
 \hline
{2.1651e-01} & $0.4$ & {5.2016e-05} & - & {2.5203e-03} & -   \\
{1.5309e-01} & $0.2\sqrt{2}$  & {1.1008e-05} & {4.481} & {1.3058e-03} & {1.897}  \\
{1.0825e-01} & $0.2$ & {2.8535e-06} & {3.896} & {6.8447e-04} & {1.864}  \\
{7.6547e-02} & $0.1\sqrt{2}$  & {6.9422e-07} & {4.079} & {3.4543e-04} & {1.973}  \\
\hline
\end{tabular}
\caption{Errors and experimental orders of convergence for Example 4 \label{tab:ex3d}}
\end{table}
\end{center}


\section{Appendix}

\begin{lemma} \label{extpde}
Suppose that $u$ is a smooth solution of (\ref{surfeq}) and denote by $u^e$ the extension defined in (\ref{ext}). Then
$u^e$ is a solution of (\ref{eq:ue2}).
\end{lemma}
{\bf Proof:}
We use the notation introduced in Section \ref{extension} and begin by  deriving a formula for $\tilde{p}(x,t)$ for  $x \in U_{\delta}(t), t \in [0,T]$. Define 
\begin{displaymath}
  \eta(\tau) := F_t(p(x,t), (1-\tau) \phi(x,t)), \quad \tau \in [0,1].
\end{displaymath}
Recalling (\ref{odesys}) and the definition of $F_t$ we have
\begin{displaymath}
  \eta'(\tau)
  = -\phi(x,t) \frac{ \nabla \phi( \gamma_{p(x,t),t} ((1-\tau) \phi(x,t)),t)}{ \abs{ \nabla \phi( \gamma_{p(x,t),t} ( ( 1 - \tau ) \phi(x,t)),t) }^2 }.
\end{displaymath}
Observing that $\gamma_{p(x,t),t} ( \phi(x,t) ) = F_t( p(x,t), \phi(x,t) ) = x$ and using similar arguments to calculate $\eta''(\tau)$ we find for
$k=1,\ldots,n+1$ that
\begin{eqnarray*}
\eta_k'(0) & = & - \phi(x,t) \frac{ \phi_{x_k}(x,t)}{ \abs{ \nabla \phi(x,t) }^2}, \\
\eta_k''(0) & = &  \phi(x,t)^2 \, \sum_{l,r=1}^{n+1} \left(
    \delta_{kr} - \frac{2 \phi_{x_k}(x,t) \phi_{x_r}(x,t)}{ \abs{ \nabla \phi(x,t) }^2} \right)
    \frac{ \phi_{x_l}(x,t) \phi_{x_l x_r}(x,t)}{\abs{ \nabla \phi(x,t) }^4}.
\end{eqnarray*}
Since $\eta(1)=F_t(p(x,t),0)=\Phi(p(x,t),t)= \tilde{p}(x,t)$, $\eta(0)=x$ we deduce with the help of Taylor's theorem that for $k=1,\ldots,n+1$
\begin{eqnarray}
\lefteqn{
\tilde{p}_k(x,t)  = x_k - \phi(x,t) \frac{\phi_{x_k}(x,t)}{| \nabla \phi(x,t) |^2} } \label{pk} \\
& & + \frac{1}{2} \phi(x,t)^2 \sum_{l,r=1}^{n+1} \Bigl( \delta_{kr}-
    \frac{2 \phi_{x_k}(x,t) \phi_{x_r}(x,t)}{| \nabla \phi(x,t) |^2} \Bigr) \frac{\phi_{x_l}(x,t) \phi_{x_l x_r}(x,t)}{| \nabla \phi(x,t) |^4}
+ \phi(x,t)^3 r_k(x,t), \nonumber
\end{eqnarray}
where $r_k$ are smooth functions. Starting from (\ref{pk}) it is not difficult to derive formulae for $\tilde{p}_{x_i}, \tilde{p}_{x_i x_j}$ (cf. (2.9), (2.10) in \cite{DER14})
and hence to deduce from (\ref{zederiv1}) and (\ref{zederiv2}) that
\begin{eqnarray}
\nabla u^e(x,t) & = & (I+\phi(x,t) A(x,t)) \nabla_{\Gamma}u( \tilde{p}(x,t),t) \label{zegrad} \\[2mm]
\lefteqn{ \hspace{-1.5cm} \frac{1}{| \nabla \phi(x,t)|} \nabla \cdot \bigl( | \nabla \phi(x,t)| \nabla u^e(x,t) \bigr)  = (\Delta_{\Gamma} u)(\tilde{p}(x,t),t) } \label{zediv} \\
& &   + \phi(x,t) \Bigl( \sum_{k,l=1}^{n+1} b_{lk}(x,t) \underline{D}_l \underline{D}_k u(\tilde{p}(x,t),t)   + \sum_{k=1}^{n+1} c_k(x) \underline{D}_k u(\tilde{p}(x,t),t) \Bigr),
  \nonumber
\end{eqnarray}
where $A=(a_{ik}),b_{lk}$ and $c_k$ are smooth. Furthermore, differentiating (\ref{pk}) with respect to $t$ we find that
\begin{displaymath}
\tilde{p}_t(x,t) = - \frac{\phi_t(x,t)}{| \nabla \phi(x,t) |^2} \nabla \phi(x,t) + \phi(x,t) q(x,t)= V(x,t) \nu(x,t) + \phi(x,t) q(x,t), \quad q \mbox{ smooth}
\end{displaymath}
we infer from (\ref{zegrad}),  (\ref{zemat}) and (\ref{tanggrad0}) that
\begin{eqnarray} 
\lefteqn{
 \partial^{\bullet}_t u^e(x,t)  =  u^e_t(x,t) + \bigl( \vel(x,t), \nabla u^e(x,t) \bigr) = u^e_t(x,t) + \bigl( \vel(x,t), (I+\phi(x,t) A(x,t)) \nabla_{\Gamma}u( \tilde{p}(x,t),t) \bigr) } \nono  \\
 & = &  \partial^{\bullet}_t u (\tilde{p}(x,t),t) + \bigl( (\vel(x,t) - \vel(\tilde{p}(x,t),t)), \nabla_{\Gamma} u(\tilde{p}(x,t),t) \bigr) \nono \\
 & & + V(x,t) \bigl( (\nu(x,t) - \nu(\tilde{p}(x,t),t)), \nabla_{\Gamma} u(\tilde{p}(x,t),t) \bigr)
  + \phi(x,t) \bigl( q(x,t) + A(x,t)^T \vel(x,t), \nabla_{\Gamma} u(\tilde{p}(x,t),t) \bigr) \nono \\
 & = & \partial^{\bullet}_t u (\tilde{p}(x,t),t) + \phi(x,t) r(x,t) \label{zet}
 \end{eqnarray}
for some smooth function $r$. Finally, since $\nabla_{\phi} \cdot \vel(\tilde{p}(x,t),t) = \nabla_{\Gamma} \cdot \vel (\tilde{p}(x,t),t)$ we have
\begin{eqnarray}
u^e(x,t) \nabla_{\phi} \cdot \vel(x,t) & = &  u(\tilde{p}(x,t),t) \nabla_{\phi} \cdot \vel(\tilde{p}(x,t),t) + u(\tilde{p}(x,t),t) \bigl( \nabla_{\phi} \cdot \vel(x,t) -
\nabla_{\phi} \cdot \vel(\tilde{p}(x,t),t) \bigr) \nono \\
& = & u(\tilde{p}(x,t),t) \nabla_{\Gamma} \cdot \vel (\tilde{p}(x,t),t) + \phi(x,t) \bar{r}(x,t).  \label{div}
\end{eqnarray}
Combining (\ref{zediv})--(\ref{div})  we deduce that the extension $u^e$ of a function $u$ solving (\ref{surfeq}) satisfies 
(\ref{eq:ue2}). \qed \\[3mm]


{\bf Acknowledgements}

The authors would like to thank the Isaac Newton Institute for
Mathematical Sciences for its hospitality during the programme 
{\it Coupling Geometric PDEs with Physics for Cell Morphology, Motility and Pattern Formation}
supported by EPSRC Grant Number EP/K032208/1 and VS gratefully acknowledges the support of the Leverhulme Trust  Research Project Grant (RPG-2014-149). 


\begin{thebibliography}{99}
\bibitem{AS03}  Adalsteinsson, D., Sethian, J.A.: Transport and diffusion of material quantities on propagating interfaces via level set methods. J. Comput. Phys. {\bf 185}, 271–-288 (2003).
\bibitem{BNS} Barrett, J.W and N\"urnberg, R. and Styles, V.: Finite element approximation of a phase field model for void electromigration. SIAM J. Numer. Anal. {\bf 46}, 738–-772 (2004).
\bibitem{BCOS} Bertalmio, M., Cheng, L.T., Osher, S., Sapiro, G.: Variational problems and partial differential equations on implicit surfaces: The framework and
examples in image processing and pattern formation. J. Comput. Phys. {\bf 174}, 759--780 (2001).
\bibitem{Bu09} Burger, M: Finite element approximation of elliptic partial differential equations on implicit surfaces. Comput. Vis. Sci. {\bf 12}, 87--100 (2009).
\bibitem{CFP97} Cahn, J.W., Fife, P., Penrose, O.: A phase field model for diffusion induced grain boundary motion. Acta Mater. {\bf 45}, 4397–-4413 (1997).
\bibitem{DDEH10} Deckelnick, K., Dziuk, G.,  Elliott, C.M.,  Heine, C.-J.:
An h-narrow band finite element method for elliptic equations on implicit surfaces.
IMA Journal of Numerical Analysis {\bf 30}, 351--376 (2010). 
\bibitem{DER14} Deckelnick, K., Elliott, C.M., Ranner, T.: Unfitted finite element methods using bulk meshes for surface partial differential equations. SIAM J. Numer. Anal. {\bf 52}, 2137–-€"2162 (2014). 
\bibitem{DE07}  Dziuk, G., Elliott, C.M.:  Finite elements on evolving surfaces. IMA J. Numer. Anal. {\bf 27}, 262–-292 (2007).
\bibitem{DE08}  Dziuk, G., Elliott, C.M.:  Eulerian finite element method for parabolic PDEs on implicit surfaces.  Interfaces Free Bound. {\bf 10}, 119-–138 (2008). 
\bibitem{DE10}  Dziuk, G., Elliott, C.M.:  An Eulerian approach to transport and diffusion on evolving implicit surfaces. Comput. Vis. Sci. {\bf 13}, 17–-28 (2010).
\bibitem{DE12}  Dziuk, G., Elliott, C.M.:  A fully discrete evolving surface finite element method.  SIAM J. Numer. Anal. {\bf 50}, 2677-–2694 (2012).
\bibitem{DE13}  Dziuk, G., Elliott, C.M.:  L$^2$–-estimates for the evolving surface finite element method. Math. Comp. {\bf 82},  1-–24 (2013).
\bibitem{DzE13}  Dziuk, G., Elliott, C.M.: Finite element methods for surface PDEs. Acta Numer. {\bf 22}, 289--396 (2013). 
\bibitem{DLM12} Dziuk, G., Lubich, C., Mansour, D.: Runge-Kutta time discretization of parabolic differential equations on evolving surfaces. IMA J. Numer. Anal. {\bf 32}, 394–-416 (2012).
\bibitem{ES09}  Elliott, C.M.,  Stinner, B.:  Analysis of a diffuse interface approach to an advection diffusion equation on a moving surface. Math. Models Methods Appl. Sci. {\bf 19},  
787–-802 (2009).
\bibitem{ESSW11}  Elliott, C.M., Stinner, B., Styles, V., Welford, R.: Numerical computation of advection and diffusion on evolving diffuse interfaces. IMA J. Numer. Anal. {\bf 31}, 
786-–812 (2011)
\bibitem{HLZ15}  Hansbo, P., Larson, M.G., Zahedi, S.: Characteristic cut finite element methods for convection–-diffusion problems on time dependent surfaces.
Comput. Methods Appl. Mech. Engrg. {\bf 293},  431-–461 (2015).
\bibitem{LNR11} Lenz, M., Nemadjieu, S.F., Rumpf, M.:  A convergent finite volume scheme for diffusion on evolving surfaces. SIAM J. Numer. Anal. {\bf 49},  15–-37 (2011).
\bibitem{LMV13} Lubich, C., Mansour, D.,  Venkataraman, C.: Backward difference time discretization of parabolic differential equations on evolving surfaces. IMA J. Numer. Anal. 
{\bf 33},  1365–-1385 (2013). 
\bibitem{ORX14}  Olshanskii, M.A.,  Reusken, A., Xu, X.:  An Eulerian space-time finite element method for diffusion problems on evolving surfaces. SIAM J. Numer. Anal. {\bf 52}, 1354–-1377 (2014).
\bibitem{OR14}   Olshanskii, M.A.,  Reusken, A.: Error analysis of a space-time finite element method for solving PDEs on evolving surfaces. SIAM J. Numer. Anal. {\bf 52},  2092–-2120 (2014).
\bibitem{RV06} R\"atz, A., Voigt, A.: PDE's on surfaces - a diffuse interface approach. Comm. Math. Sci. {\bf 4}, 575--590 (2006).
\bibitem{alberta}  Schmidt, A. and Siebert, K.G.: 
Design of adaptive finite element software. The finite element toolbox ALBERTA.  
Lecture Notes in Computational Science and Engineering 42,
Springer-Verlag, Berlin, (2005). 
\bibitem{S90}  Stone, H.A.:  A simple derivation of the time–-dependent convective–-diffusive equation for surfactant transport along a deforming interface. Phys. Fluids A {\bf 2}, 111–-112 (1990).
\bibitem{TLLWV09}  Teigen, K.E.,  Li, X.,  Lowengrub, J., Wang, F., Voigt, A.:  A diffuse–-interface approach for modeling transport, diffusion and adsorption/desorption of material quantities on a deformable interface. Commun. Math. Sci. {\bf 7},  1009–-1037 (2009).
\bibitem{XZ03} Xu, J.J., Zhao, H.K.:  An Eulerian formulation for solving partial differential equations along a moving interface. J. Sci. Comput. {\bf 19}, 573–-594 (2003).
\end{thebibliography}
\end{document}